\RequirePackage{ifpdf}
\ifpdf 
\documentclass[pdftex]{sigma}
\else
\documentclass{sigma}
\fi

\numberwithin{equation}{section}

\newtheorem{Theorem}{Theorem}[section]
\newtheorem{Corollary}[Theorem]{Corollary}
\newtheorem{Lemma}[Theorem]{Lemma}
\newtheorem{Proposition}[Theorem]{Proposition}
{ \theoremstyle{definition}
\newtheorem{Definition}[Theorem]{Definition}

\newtheorem{Remark}[Theorem]{Remark} }

\usepackage{upgreek}
\usepackage[mathscr]{euscript}
\usepackage{cleveref}

\begin{document}

\newcommand{\arXivNumber}{1912.06067}

\renewcommand{\PaperNumber}{021}

\FirstPageHeading

\ShortArticleName{Parameter Permutation Symmetry in Particle Systems and Random Polymers}

\ArticleName{Parameter Permutation Symmetry\\ in Particle Systems and Random Polymers}

\Author{Leonid PETROV~$^{{\rm ab}}$}

\AuthorNameForHeading{L.~Petrov}

\Address{$^{\rm a)}$~University of Virginia, Department of Mathematics,\\
\hphantom{$^{\rm a)}$}~141 Cabell Drive, Kerchof Hall, P.O.~Box 400137, Charlottesville, VA 22904, USA}
\EmailD{\href{mailto:lenia.petrov@gmail.com}{lenia.petrov@gmail.com}}
\URLaddressD{\url{http://lpetrov.cc}}
\Address{$^{\rm b)}$~Institute for Information Transmission Problems,\\
\hphantom{$^{\rm b)}$}~Bolshoy Karetny per.~19, Moscow, 127994, Russia}

\ArticleDates{Received October 26, 2020, in final form February 20, 2021; Published online March 06, 2021}

\Abstract{Many integrable stochastic particle systems in one space dimension (such as TASEP~-- totally asymmetric simple exclusion process~-- and its various deformations, with a notable exception of ASEP) remain integrable when we equip each particle $x_i$ with its own jump rate parameter $\nu_i$. It is a consequence of integrability that the distribution of~each particle $x_n(t)$ in a system started from the step initial configuration depends on the parameters $\nu_j$, $j\le n$, in a symmetric way. A transposition $\nu_n \leftrightarrow \nu_{n+1}$ of the para\-me\-ters thus affects only the distribution of $x_n(t)$. For $q$-Hahn TASEP and its degenerations (\mbox{$q$-TASEP} and directed beta polymer) we realize the transposition $\nu_n \leftrightarrow \nu_{n+1}$ as an explicit Markov swap operator acting on the single particle~$x_n(t)$. For beta polymer, the swap operator can be interpreted as a simple modification of the lattice on which the polymer is considered. Our main tools are Markov duality and contour integral formulas for joint moments. In~particular, our constructions lead to a continuous time Markov process $\mathsf{Q}^{(\mathsf{t})}$ preserving the time $\mathsf{t}$ distribution of the $q$-TASEP (with step initial configuration, where $\mathsf{t}\in \mathbb{R}_{>0}$ is fixed). The dual system is a certain transient modification of the stochastic $q$-Boson system. We~identify asymptotic survival probabilities of this transient process with $q$-moments of~the~$q$-TASEP, and use this to show the convergence of the process $\mathsf{Q}^{(\mathsf{t})}$ with arbitrary initial data to its stationary distribution. Setting $q=0$, we recover the results about the usual TASEP established recently in~[arXiv:1907.09155] by a different approach based on Gibbs ensembles of interlacing particles in two dimensions.}

\Keywords{$q$-TASEP; stochastic $q$-Boson system; stationary distribution; coordinate Bethe ansatz; $q$-Hahn TASEP}

\Classification{82C22; 60C05; 60J27}

\section{Introduction}

\subsection{Overview}

In the past two decades,
integrable stochastic interacting particle
systems in one space dimension
have been crucial in explicitly describing
new universal asymptotic phenomena, most notably those corresponding
to the Kardar--Parisi--Zhang universality class
\cite{CorwinKPZ,Corwin2016Notices, halpin2015kpzCocktail,QuastelSpohnKPZ2015}.
By \emph{integrability} in a stochastic system we mean
the presence of exact formulas for probability distributions
for a wide class of observables.
Asymptotic (long time and large space)
behavior of the system can be recovered by an analysis of these formulas.
Initial successes with integrable stochastic particle systems
were achieved through the use of determinantal point process
techniques, e.g., see~\cite{johansson2000shape} for the asymptotic fluctuations of
TASEP (totally asymmetric simple exclusion process).
More recently new tools borrowed from quantum integrability, Bethe ansatz, and/or symmetric functions
were applied to deformations of TASEP and related models:
\begin{itemize}\itemsep=0pt
	\item
		ASEP, in which particles can jump in both directions,
		but with different rates\footnote{We say that an event
		in continuous time happens at rate $\alpha$ if
		$\mathbb{P}(\textnormal{waiting time till the event occurs}>t)={\rm e}^{-\alpha t}$ for all
		$t\in \mathbb{R}_{\ge0}$.}~\cite{TW_ASEP1, TW_ASEP2};
	\item
		random polymers such as the semi-discrete directed Brownian polymer~\cite{Oconnell2009_Toda},
		log-gamma polymer~\cite{COSZ2011,OSZ2012, Seppalainen2012},
		or beta type polymers		\cite{CorwinBarraquand2015Beta,BufetovMucciconiPetrov2018,CMP_qHahn_Push, MucciconiPetrov2020, thieryLD2015integrable};
	\item
		$q$-TASEP and $q$-Hahn TASEP,
		in which particles
		jump in one direction, but with $q$-defor\-med jump rates
		\cite{BorodinCorwin2011Macdonald,BorodinCorwinSasamoto2012,Corwin2014qmunu, FerrariVeto2013,Povolotsky2013}.
\end{itemize}

All these and several other integrable models
can be unified under the umbrella of
stochastic vertex models
\cite{BCG6V, BorodinPetrov2016_Hom_Lectures, borodin_wheeler2018coloured,CorwinPetrov2015}.

Ever since the original works on TASEP
around the year 2000
it was clear
\cite{Gravner-Tracy-Widom-2002a,Its-Tracy-Widom-2001a}
that
integrability of some particle systems like TASEP
is preserved in the presence of countably many extra
parameters, for example, when each particle
is equipped with its own jump rate.
We will refer to such more general systems as \emph{multiparameter} ones.
This notion should be contrasted with the $q$-deformation by means of just one extra parameter
which takes TASEP
to $q$-TASEP. The latter is much more subtle and relies on passing to a
deformed algebraic structure~-- for the $q$-TASEP, one replaces the Schur symmetric functions
with the $q$-Whittaker ones.

It should be noted that
TASEP in inhomogeneous space (when the jump rate of a particle depends on its location)
does not seem to be integrable~\cite{costin2012blockage, janowsky1992slow_bond, seppalainen2001slow_bond} (cf.\ recent asymptotic fluctuation results
\cite{Basuetal2014_slowbond,basu2017invariant} requiring very delicate asymptotic analysis).
Moreover, it is not known whether ASEP has any
integrable multiparameter
deformations.
The stochastic six vertex model~\cite{BCG6V,GwaSpohn1992}
scales to ASEP and admits such a multiparameter deformation
\cite{BorodinPetrov2016inhom}, but this deformation is destroyed by the scaling.
Recently other families of
spatially inhomogeneous
integrable stochastic particle systems in one and two space dimensions
were
studied in~\cite{theodoros2019_determ, BorodinPetrov2016Exp,SaenzKnizelPetrov2018,Petrov2017push}.

All known multiparameter integrable stochastic particle systems
display a common feature. Namely, certain joint distributions in these
systems are symmetric under (suitably restricted classes of) permutations of the parameters.
This symmetry is far from being evident from the beginning, and
is often observed only as a consequence of explicit formulas.
The main goal of~the present paper is to \textit{explore
probabilistic consequences of
parameter symmetries in integrable particle systems}.

Recently a number of other
papers investigating symmetries of multiparameter
integrable
stochastic particle systems and vertex models have appeared
\cite{borodin2019shift,corwin2020invariance,dauvergne2020hidden,galashin2020symmetries}.
So far it is not clear whether those results
have any direct connection to the results of
the present paper.

\subsection[Distributional symmetry of the q-Hahn TASEP]
{Distributional symmetry of the $\boldsymbol q$-Hahn TASEP}
The most general system
we consider is the $q$-Hahn TASEP
started from the
\emph{step initial con\-fi\-gu\-ra\-tion}
$x_n(0)=-n$, $n=1,2,\ldots $.
That is, every site of $\mathbb{Z}_{<0}$ is occupied by a particle,
and every site of $\mathbb{Z}_{\ge0}$ is empty.
Throughout
the paper
we denote this configuration by
$\mathsf{step}$ for short.

\begin{figure}[htbp]\centering
	\includegraphics{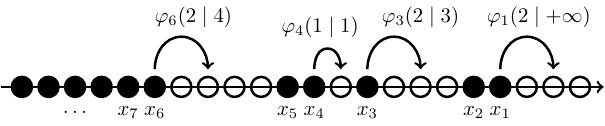}
	\caption{An example of a one-step transition in the $q$-Hahn TASEP,
		together with the corresponding probabilities for each particle.
		Here $\varphi_i\equiv \varphi_{q,\gamma\nu_i,\nu_i}$.} 	\label{fig:qhahn_tasep}
\end{figure}

The $q$-Hahn TASEP was
introduced in~\cite{Povolotsky2013}
and studied in~\cite{BCPS2014_arXiv_v4, Corwin2014qmunu,Veto2014qhahn}.
Its multiparameter deformation
appears in~\cite{BorodinPetrov2016inhom}.
Under this deformation, each particle $x_n$
carries its own parameter
$\nu_n\in(0,1)$ which determines the jump distribution of the particle.
The $q$-Hahn TASEP is a~discrete time Markov
process on particle configurations in $\mathbb{Z}$.
At each time step, every particle~$x_i$ independently
jumps to the right by $j$ steps with probability
\begin{gather*}
	\varphi_{q,\gamma\nu_i,\nu_i}(j\,|\, x_{i-1}-x_i-1),
	\qquad
	j\in \left\{ 0,1,\ldots,x_{i-1}-x_i-1 \right\},
\end{gather*}
where $x_0=+\infty$, by agreement.
Here $\varphi$ is the $q$-deformed beta-binomial distribution
(Definition~\ref{def:phi_distribution}).
See Fig.~\ref{fig:qhahn_tasep} for an illustration.

\begin{figure}[htpb]
	\centering
	\includegraphics{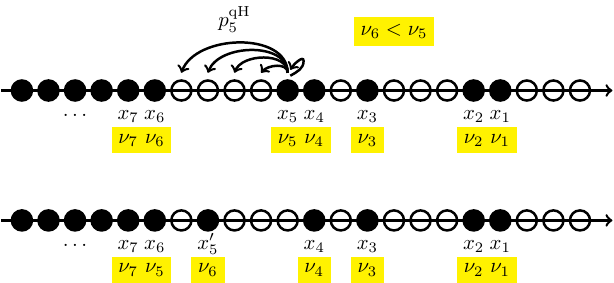}
	\caption{An example of the swap operator $p_5^{\mathrm{qH}}$ acting on the fifth particle
		in the $q$-Hahn TASEP (at~an~arbitrary time $t\in \mathbb{Z}_{\ge0}$).
		Arrows show possible new locations of~$x_5$ (note that
		with some probability it can stay in the same location).
		The resulting configuration (below) is distributed as the $q$-Hahn TASEP at the same time,
		but with the swapped parameters $\nu_5 \leftrightarrow \nu_6$.
		The distributional identity holds only if~$\nu_6<\nu_5$ before the swap.}
	\label{fig:intro_fig_qhahn}
\end{figure}

The distribution of each particle $x_n(t)$ at any time moment
$t\in \mathbb{Z}_{\ge0}$
in the $q$-Hahn TASEP started from $\mathsf{step}$
depends on the parameters $\nu_1,\ldots,\nu_n $ in a symmetric way.
We check this~sym\-metry using exact formulas in Section~\ref{sub:qhahn_defn}.
The main structural result of the present paper~is
\begin{Theorem}[Theorem~\ref{thm:qhahn_swap} in the text]
	\label{thm:intro_qhahn_swap}
	The elementary transposition
	$\nu_n \leftrightarrow \nu_{n+1}$, \mbox{$\nu_{n+1}<\nu_n$},
	of two neighboring parameters
	in the $q$-Hahn TASEP started from $\mathsf{step}$ is
	equivalent in distribution
	to the action of an explicit Markov swap operator~$p_n^{\mathrm{qH}}$
	on the particle $x_n$.
	This operator moves $x_n$ to a random new location $x_n'$ chosen
	with probability
	\begin{gather*}
		\varphi_{q,\frac{\nu_{n+1}}{\nu_n},\nu_{n+1}}(x_n'-x_{n+1}-1\,|\, x_n-x_{n+1}-1),\qquad
		x_n'\in \{x_{n+1}+1,\ldots,x_n-1,x_n \},
	\end{gather*}
	where
	$\varphi$ is the $q$-deformed beta-binomial distribution
	$($Definition~$\ref{def:phi_distribution})$.
	The equivalence in distribution holds at any fixed time
	$t \in \mathbb{Z}_{\ge0}$ in the $q$-Hahn TASEP,
	while the swap operator $p_n^{\mathrm{qH}}$ does not depend on~$t$.
	See Fig.~$\ref{fig:intro_fig_qhahn}$ for an illustration.
\end{Theorem}

We prove this result
in Section~\ref{sub:cond_distr_qhahn}
using $q$-moment contour integral formulas
and duality results\footnote{Sometimes,
especially in the context of random polymers,
this set of tools is referred to as ``rigorous replica method''.}
for the $q$-Hahn TASEP.
Let us make a couple of remarks on the generality of the result and our methods.

First, note that there are certain other classes of
initial data
(for example, half-stationary)
for which the $q$-Hahn TASEP displays parameter symmetry.
Moreover, via the spectral theory of
\cite{BCPS2014_arXiv_v4}
one sees that for fairly general initial
data the swap operators simultaneously applied
to the particle system and the initial distribution
lead to permutations of the parameters as in Theorem~\ref{thm:intro_qhahn_swap}.
For simplicity, in this paper we focus only on the
step initial configuration.

Second, we expect that under suitable modifications
the results of the present paper could carry over to the
stochastic six vertex model~\cite{BCG6V}
and the higher spin stochastic six vertex model~\cite{BorodinPetrov2016inhom, CorwinPetrov2015}.
However, as duality relations for these vertex models
are more involved than the ones for the $q$-Hahn TASEP,
it is not immediately clear how to extend the methods of the present
paper to the vertex models. Therefore, we restrict
out attention here to the $q$-Hahn TASEP.

\subsection{Applications}

We explore a number of interesting consequences of the
distributional symmetry of the $q$-Hahn TASEP
realized by the swap operators.
Let us briefly describe them.

We take a \emph{continuous time limit} of the
$q$-Hahn TASEP and the swap operators.
Denote by~$\mathscr{M}^{\mathrm{qH}}_{q,\nu;\mathsf{t}}$
the distribution of the parameter
homogeneous
(i.e., $\nu_n\equiv \nu$),
continuous time \mbox{$q$-Hahn} TASEP
at time $\mathsf{t}\in \mathbb{R}_{\ge0}$ started from $\mathsf{step}$
(see Section~\ref{sub:cont_time_qhahn} for a detailed definition).
The~$q$-Hahn TASEP
evolution acts on
$\mathscr{M}^{\mathrm{qH}}_{q,\nu;\mathsf{t}}$
by increasing the time parameter $\mathsf{t}$.
We find that a~suitable continuous limit as $r\to1$ of
the swap operators with $\nu_n=\nu r^{n-1}$
produces a (time-inhomogeneous) continuous time Markov process
$\mathcal{B}^{\mathrm{qH}}$ on particle configurations.
Starting from a~random particle configuration
distributed as
$\mathscr{M}^{\mathrm{qH}}_{q,\nu;\mathsf{t}}$
and running the process
$\mathcal{B}^{\mathrm{qH}}$
for time $\tau\ge0$, we get a~con\-fi\-gu\-ra\-tion distributed as
$\mathscr{M}^{\mathrm{qH}}_{q,\nu {\rm e}^{-\tau};\mathsf{t} {\rm e}^{-\tau}}$,
that is, in which both parameters $\nu$ and $\mathsf{t}$ are rescaled.
See~Theorem~\ref{thm:action_on_qhahn_distributions}
for a detailed formulation
and
Fig.~\ref{fig:two_times_qhahn_diagram} for an illustration of the two actions.
When $\nu=0$, the backward process
becomes time-homogeneous, and we discuss this
case in more detail in the next
Section~\ref{sub:intro_stat_tasep_briefly}.

When $q=\nu=0$, Theorem~\ref{thm:action_on_qhahn_distributions} recovers one of the main results of
the recent work~\cite{PetrovSaenz2019backTASEP}
on~the existence of a
time-homogeneous,
continuous time process mapping the distributions
of the usual TASEP
back in time. We remark that the proof of
this result following from the present paper
is completely different
from the argument given in~\cite{PetrovSaenz2019backTASEP}.
The latter went through
the well-known connection of the TASEP distribution
and a~Schur process~\cite{okounkov2003correlation}
on interlacing arrays
(about this connection see, e.g.,~\cite{BorFerr2008DF}).
For Schur processes, the two-dimensional version of the
swap operator
is accessible by elementary means.

In a scaling limit $q$, $\nu_n\to1$, the $q$-Hahn TASEP turns into the beta
polymer model introduced in~\cite{CorwinBarraquand2015Beta}.
In Section~\ref{sec:beta_polymer} we construct swap operators
for the multiparameter version of the beta polymer model.
The argument is formally independent of the rest of the paper,
but proceeds through the same steps.
For polymers, the swap operator can be realized
as a simple modification of the lattice on which the beta polymer is defined.
See Theorem~\ref{thm:beta_polymer_back}
for a detailed formulation of the result, and
Figs.~\ref{fig:beta_modified_lattice} and~\ref{fig:beta_modified_lattice2}
for illustrations of lattice modifications.

\subsection[Stationary dynamics on the q-TASEP distribution]{Stationary dynamics on the $\boldsymbol q$-TASEP distribution}\label{sub:intro_stat_tasep_briefly}

The last application concerns $q$-TASEP
\cite{BorodinCorwin2011Macdonald,BorodinCorwinSasamoto2012}, which is a
$\nu=0$ degeneration of the $q$-Hahn TASEP.
Let us focus on this case in more detail.
Under the $q$-TASEP, each particle $x_n$ jumps to~the right by one in continuous time
at rate $1-q^{x_{n-1}-x_{n}-1}$, where $x_0=+\infty$, by agreement.
In~particular, we take the homogeneous $q$-TASEP in which all particles behave in the same manner.
Denote by
$\mathscr{M}^{\mathrm{qT}}_{q;\mathsf{t}}$
the time $\mathsf{t}$
distribution of this
continuous time $q$-TASEP
started from the step initial configuration $\mathsf{step}$.

When $\nu=0$, Theorem~\ref{thm:action_on_qhahn_distributions} produces a
new time-homogeneous, continuous time Markov process
which we denote by $\mathsf{Q}^{(\mathsf{t})}$, with the
following properties:
\begin{itemize}\itemsep=0pt
	\item
		The process
		$\mathsf{Q}^{(\mathsf{t})}$
		is a combination of two independent dynamics:
		the $q$-TASEP evolution, and the
		(slowed down by the factor of $\mathsf{t}$)
		backward $q$-TASEP evolution. The latter is a~suitable degeneration of the
		backward $q$-Hahn process $\mathcal{B}^{\mathrm{qH}}$. Under this degeneration,
		the backward process becomes time-homogeneous.
		See Section~\ref{sub:stat_qtasep_markov_gen} for the full definition
		of the process $\mathsf{Q}^{(\mathsf{t})}$.
	\item
		(Proposition~\ref{prop:stationary_qtasep_process})
		The process $\mathsf{Q}^{(\mathsf{t})}$ preserves the
		distribution
		$\mathscr{M}^{\mathrm{qT}}_{q;\mathsf{t}}$.
		Here the time parameter $\mathsf{t}\in \mathbb{R}_{\ge0}$
		of the $q$-TASEP distribution
		$\mathscr{M}^{\mathrm{qT}}_{q;\mathsf{t}}$
		is fixed and is incorporated
		into the definition of the stationary process
		$\mathsf{Q}^{(\mathsf{t})}$.
	\item
		(Theorem~\ref{thm:qtasep_converges_to_stationary})
		Start the process
		$\mathsf{Q}^{(\mathsf{t})}$ from an arbitrary
		particle configuration on $\mathbb{Z}$
		which is empty far to the
		right, densely packed far to the left, and is balanced
		(in the sense that the number of holes to the left
		of zero equals the number of particles to the right of
		zero).
		Then in the long time limit the distribution of this
		process converges to the $q$-TASEP distribution
		$\mathscr{M}^{\mathrm{qT}}_{q;\mathsf{t}}$.
\end{itemize}

We establish Theorem~\ref{thm:qtasep_converges_to_stationary}
by making use of duality for the stationary process
$\mathsf{Q}^{(\mathsf{t})}$ which extends the duality
between the
$q$-TASEP
and the stochastic $q$-Boson process
from~\cite{BorodinCorwinSasamoto2012}
(the stochastic $q$-Boson process dates back to~\cite{SasamotoWadati1998}).
In fact, we are able to use the same duality functional (corresponding to joint $q$-moments)
for $\mathsf{Q}^{(\mathsf{t})}$.
As a result we find that the process dual to
$\mathsf{Q}^{(\mathsf{t})}$
is a new transient modification of the stochastic $q$-Boson process.
The long time limit of this transient process
is readily accessible, and
Theorem~\ref{thm:qtasep_converges_to_stationary}
follows by matching the
long time behavior of all $q$-moments of the stationary dynamics
$\mathsf{Q}^{(\mathsf{t})}$ (with an arbitrary initial configuration)
with those of the $q$-TASEP (with the step initial configuration).

Let us illustrate the transient modification in the simplest case of
the first $q$-moment.
Consider the continuous time random walk $n(\mathsf{t})$ on $\mathbb{Z}_{\ge0}$
which jumps from $k$ to $k-1$, $k\ge1$, at rate $1-q$.
When the walk reaches zero, it stops.
From the $q$-TASEP duality
\cite{BorodinCorwinSasamoto2012} we have
\begin{gather}
	\label{eq:qTASEP_single_particle_intro}
	\mathbb{E}^{\mathrm{qT}}_{\mathsf{step}}
	\,q^{x_n(\mathsf{t})+n}=
	\mathbb{P}(n(\mathsf{t})>0\,|\, n(0)=n),\qquad
	n=1,2,\ldots .
\end{gather}
Here the left-hand side is the expectation over the $q$-TASEP distribution
$\mathscr{M}^{\mathrm{qT}}_{q;\mathsf{t}}$, and the
right-hand side corresponds to the random walk $n(\mathsf{t})$.
Similarly to~\eqref{eq:qTASEP_single_particle_intro},
joint $q$-moments of the $q$-TASEP are governed by a
multiparticle version of the process $n(\mathsf{t})$
--- the stochastic $q$-Boson system.

Let us now fix the $q$-TASEP time parameter $\mathsf{t}\in \mathbb{R}_{>0}$,
and consider a different random walk $n^{(\mathsf{t})}(\tau)$ on $\mathbb{Z}_{\ge0}$
(here $\tau$ is the new continuous time variable)
with the following jump rates:
\begin{gather*}
	\mathrm{rate}(k\to k-1)=1-q,
	\qquad
	\mathrm{rate}(k-1\to k)=\frac{k-1}{\mathsf{t}},
	\qquad
	k=1,2,\ldots .
\end{gather*}
This process has a single absorbing state $0$,
and otherwise is transient. In other words, after a large time $\tau$, the particle
$n^{(\mathsf{t})}(\tau)$ is either at $0$, or runs off to infinity.
Note however that this process does not make
infinitely many jumps in finite time.
The duality for the stationary process $\mathsf{Q}^{(\mathsf{t})}$
which we prove in this paper
states (in the simplest case)
that
\begin{gather}
	\label{eq:stat_qTASEP_single_particle_intro}
	\mathbb{E}^{\mathrm{stat}(\mathsf{t})}_{\mathsf{step}}
	q^{x_n(\tau)+n}=
	\mathbb{P}\big(n^{(\mathsf{t})}(\tau)>0\,|\, n^{(\mathsf{t})}(0)=n\big),
	\qquad n=1,2,\ldots.
\end{gather}
Here the left-hand side is the expectation over the
stationary process started from $\mathsf{step}$, and
the right-hand side
may be called the \emph{survival probability}
(up to time $\tau$) of
the transient random walk $n^{(\mathsf{t})}$.
See
Corollary~\ref{cor:duality_Q_expectations}
for the general statement which connects
joint $q$-moments of the stationary process
$\mathsf{Q}^{(\mathsf{t})}$
with a multiparticle version of
$n^{(\mathsf{t})}(\tau)$. We call this multiparticle process the
\emph{transient stochastic $q$-Boson system}.

Taking the long time limit of~\eqref{eq:stat_qTASEP_single_particle_intro},
we see that
\begin{gather}
	\label{eq:harmonic_function_intro}
	\lim_{\tau\to+\infty}
	\mathbb{E}^{\mathrm{stat}(\mathsf{t})}_{\mathsf{step}}
	q^{x_n(\tau)+n}=
	\mathbb{P}\big(\lim_{\tau\to+\infty}
	n^{(\mathsf{t})}(\tau)=+\infty\,|\, n^{(\mathsf{t})}(0)=n\big),\qquad n=1,2,\ldots.
\end{gather}
The right-hand side is
the \emph{asymptotic survival probability} that the transient random walk
eventually runs off to infinity and is not absorbed at zero.
This probability, viewed as a function of the initial location $n$,
is a harmonic function\footnote{A harmonic function for a continuous time Markov
	process on a discrete space
	is a function which is eliminated by
	the infinitesimal generator of the process.}
for the transient random walk $n^{(\mathsf{t})}(\tau)$,
which, moreover, takes value $1$ at $n=+\infty$. This identifies
the harmonic function uniquely.
From the stationarity of $\mathscr{M}^{\mathrm{qT}}_{q;\mathsf{t}}$
under the process $\mathsf{Q}^{(\mathsf{t})}$ one can check that
\eqref{eq:qTASEP_single_particle_intro} satisfies the same harmonicity condition,
which implies that~\eqref{eq:harmonic_function_intro} equals~\eqref{eq:qTASEP_single_particle_intro}.

{\sloppy
A general multiparticle argument involves identifying
the ``correct'' harmonic function (asymptotic survival probability)
of the transient stochastic $q$-Boson system with the joint \mbox{$q$-moments} of~the~$q$-TASEP.
This identification requires additional technical steps since the space
of harmonic functions for the multiparticle process is higher-dimensional.
Along this route we obtain
the proof
that $\mathsf{Q}^{(\mathsf{t})}$
converges to its stationary distribution $\mathscr{M}^{\mathrm{qT}}_{q;\mathsf{t}}$
(Theorem~\ref{thm:qtasep_converges_to_stationary}).

}

\subsection{Outline}

In Section~\ref{sec:general_formalism} we give
general definitions related to
parameter-symmetric
particle systems
and swap operators.
In Section~\ref{sec:qhahn_TASEP} for the $q$-Hahn TASEP we present an explicit realization of the
parameter transposition in terms of a Markov swap operator
corresponding to a random jump of a single particle.
In Section~\ref{sec:qhahn_cont_time} we pass to the continuous time in the $q$-Hahn TASEP,
and obtain the $q$-Hahn backward process. This also implies the
results about the TASEP from~\cite{PetrovSaenz2019backTASEP}.
In Section~\ref{sec:stat_qtasep} we define and study the
dynamics preserving the $q$-TASEP distribution, and
show its convergence to stationarity.
In Section~\ref{sec:beta_polymer} we obtain
swap operators for the beta polymer.

\section{From symmetry to swap operators}
\label{sec:general_formalism}

This section contains an abstract discussion of
stochastic particle systems on $\mathbb{Z}$ which
depend symmetrically on their parameters.
The main notions which we use in other parts of the
paper are \emph{parameter-symmetric stochastic particle system}
and \emph{swap operators}.

\subsection{Parameter-symmetric particle systems}

Let $\mathrm{Conf}_{\rm fin}(\mathbb{Z})$ be the space of
particle configurations $\mathbf{x}=(\cdots <x_3<x_2<x_1)$,
$x_i\in \mathbb{Z}$,
which can be obtained from the \emph{step configuration}
$\mathsf{step}:=(\ldots,-3,-2,-1 )$
by finitely many operations of moving a particle to the right by one
into the nearby empty spot. The space
$\mathrm{Conf}_{\rm fin}(\mathbb{Z})$ is~countable.

By a multiparameter interacting particle system $\mathbf{x}(t)$
we mean a Markov process on $\mathrm{Conf}_{\rm fin}(\mathbb{Z})$ evolving in continuous
or discrete time,
such that $\mathbf{x}(0)=\mathsf{step}$.
Assume that this Markov process depends on countably many parameters
${\boldsymbol \nu}=\{\nu_i\}_{i\in \mathbb{Z}_{\ge1}}$.
The parameters $\nu_i$ in our situation are real, though without loss of generality
they may belong to an abstract space.
One should think that $\nu_i$ is attached to the particle $x_i$, but the distribution
of each $x_j(t)$ may depend on all of~the~$\nu_i$'s.
We denote the process depending on ${\boldsymbol \nu}$ by $\mathbf{x}^{{\boldsymbol \nu}}(t)$.
In this section we assume that all the parameters $\nu_i$ are pairwise distinct.

The infinite symmetric group
$S(\infty)=\bigcup_{n=1}^{\infty}S(n)$
acts on the parameters ${\boldsymbol \nu}$ by permutations,
$\sigma\colon {\boldsymbol \nu}\mapsto \sigma{\boldsymbol \nu}$.
Here $S(n)$ is the symmetric group which permutes only the
first parame\-ters~$\nu_1,\ldots,\nu_n$.
Let us denote by $S_n(\infty)\subset S(\infty)$
the subgroup which permutes $\nu_{n+1},\nu_{n+2},\ldots $
and maps each $\nu_i$, $1\le i\le n$, into itself. Note that
$S(n)\cap S_n(\infty)=S(n)\cap S_{n-1}(\infty)=\{e\}$,
and~$S(n+1)\cap S_{n-1}(\infty)=\{e,s_{n}\}$,
where $e$ is the identity permutation, and
$s_{n}=(n,n+1)$ is the transposition
$n\leftrightarrow n+1$.

By imposing a specific distributional symmetry
of $\mathbf{x}^{{\boldsymbol \nu}}(t)$
under the action of $S(\infty)$
on the parameters ${\boldsymbol \nu}$, we arrive at the following definition:

\begin{Definition}
	\label{def:symm_IPS}
	A multiparameter particle system $\mathbf{x}^{{\boldsymbol \nu}}(t)$ is called \emph{parameter-symmetric},
	if for all $n$ and $t$ we have
	the following equality of joint distributions:
\begin{gather}
\big(\ldots,x_{n+2}^{{\boldsymbol \nu}}(t),x_{n+1}^{{\boldsymbol \nu}}(t),
x_{n-1}^{{\boldsymbol \nu}}(t),\ldots, x_1^{{\boldsymbol \nu}}(t)\big)\nonumber
\\ \qquad
{}\stackrel{d}{=}\big(\ldots,x_{n+2}^{s_n{\boldsymbol \nu}}(t),x_{n+1}^{s_n{\boldsymbol \nu}}(t),
x_{n-1}^{s_n{\boldsymbol \nu}}(t),\ldots,x_1^{s_n{\boldsymbol \nu}}(t)\big).
\label{eq:symm_IPS}
\end{gather}
	That is, the transposition $s_n$ preserves the joint
	distribution of all particles except $x_n$.
\end{Definition}
Here is a straightforward corollary of this definition:
\begin{Corollary}
	In a parameter-symmetric particle system,
	for any $t$ and any $\sigma\in S(n)\cup S_{n}(\infty)$, the
	random variables
	$x_n^{{\boldsymbol \nu}}(t)$ and
	$x_n^{\sigma {\boldsymbol \nu}}(t)$
	have the same distribution.
\end{Corollary}

\begin{Remark}
	In Section~\ref{sec:beta_polymer} below we consider
	the beta polymer model, which may also be viewed as a
	particle system, but the particles live in $(0,1]$.
	For concreteness, in the general discussion in
	this section we stick to particle systems in $\mathbb{Z}$.
\end{Remark}

\subsection{Coupling}

Let $m_1$, $m_2$ be two probability measures
on the same measurable space $(\mathsf{E},\mathscr{F})$. A
\emph{coupling} between $m_1$ and $m_2$ is,
by definition, a measure $M=M({\rm d}z,{\rm d}z')$ on
$(\mathsf{E}\times\mathsf{E},\mathscr{F}\otimes \mathscr{F})$ whose marginals are
$m_1({\rm d}z)$ and $m_2({\rm d}z')$, respectively:
\begin{gather*}
	\int_{z'\in\mathsf{E}} M(\cdot,{\rm d}z')=m_1(\cdot),\qquad
\int_{z\in\mathsf{E}} M({\rm d}z,\cdot)=m_2(\cdot).
\end{gather*}
A coupling is not defined uniquely, but always exists
(the product measure $M=m_1\otimes m_2$ is an~example).

In the notation of the previous section, start from a parameter-symmetric
particle system~$\mathbf{x}^{{\boldsymbol \nu}}(t)$. Fix time $t$ and index $n\in \mathbb{Z}_{\ge1}$, and
consider two
distributions $\mathbf{x}^{{\boldsymbol \nu}}(t)$ and $\mathbf{x}^{s_n{\boldsymbol \nu}}(t)$
on the same countable space $\mathrm{Conf}_{\rm fin}(\mathbb{Z})$.
We would like to find a coupling $M=M_n$ between the distributions of
$\mathbf{x}^{{\boldsymbol \nu}}(t)$ and $\mathbf{x}^{s_n{\boldsymbol \nu}}(t)$ which satisfies
an additional
constraint corresponding to~\eqref{eq:symm_IPS}:
\begin{gather}
	\label{eq:symm_IPS_coupling}
	M_n\left( x_k^{{\boldsymbol \nu}}(t)=x_k^{s_n{\boldsymbol \nu}}(t)\textnormal{ for all $k\in \mathbb{Z}_{\ge1}$,
	$k\ne n$} \right)=1.
\end{gather}
Such a coupling also might not be defined uniquely.
An example of a coupling satisfying~\eqref{eq:symm_IPS_coupling}
can be obtained by adapting the basic product measure example.
For any particle configuration $\mathbf{y}=(y_1,y_2,\ldots )$ denote
$\mathbf{y}_{\hat n}:=\{y_k\colon k\ne n\}$.
Define
\begin{gather}
M_n^{\mathrm{indep}}\big(\mathbf{x}^{{\boldsymbol \nu}}(t)=\mathbf{y},\mathbf{x}^{s_n{\boldsymbol \nu}}(t)=\mathbf{z}\big)\nonumber
\\ \qquad
:=\delta(\mathbf{y}_{\hat n}=\mathbf{z}_{\hat n})
P\big(\mathbf{x}_{\hat n}^{{\boldsymbol \nu}}(t)=\mathbf{y}_{\hat n}\big)\,
P\big(x_n^{{\boldsymbol \nu}}(t)=y_n\,|\, \mathbf{x}_{\hat n}^{{\boldsymbol \nu}}(t)\big)\,
P\big(x_n^{s_n{\boldsymbol \nu}}(t)=z_n\,|\, \mathbf{x}_{\hat n}^{s_n{\boldsymbol \nu}}(t)\big).
\label{eq:symm_IPS_indep_coupling}
\end{gather}
Here $\delta(\cdot)$ is the Dirac delta,
and the two quantities
$P(\cdot\,|\, \cdot)$ are the conditional distributions
of $x_n^{{\boldsymbol \nu}}(t)$ (resp. $x_n^{s_n{\boldsymbol \nu}}(t)$)
given the locations of all other particles.
Note that
both conditional distributions
$P(\cdot\,|\, \cdot)$
in~\eqref{eq:symm_IPS_indep_coupling} are
supported on
the same interval
\begin{gather}
	\label{eq:IPS_interval_for_xn}
	I_n:=
	\big\{ x_{n+1}^{{\boldsymbol \nu}}(t)+1,x_{n+1}^{{\boldsymbol \nu}}(t)+2,\ldots,x_{n-1}^{{\boldsymbol \nu}}(t)-1 \big\}
	\subset\mathbb{Z}
\end{gather}
(if $n=1$, then, by agreement, $x_0\equiv+\infty$ and the interval is infinite; for $n\ge2$
the interval is finite).
The next statement follows from the above definitions:
\begin{Lemma}
	The distribution $M_n^{\mathrm{indep}}$ defined by~\eqref{eq:symm_IPS_indep_coupling}
	is a coupling between the distributions of~$\mathbf{x}^{{\boldsymbol \nu}}(t)$ and $\mathbf{x}^{s_n{\boldsymbol \nu}}(t)$,
	and satisfies~\eqref{eq:symm_IPS_coupling}.
\end{Lemma}

\subsection{Swap operators}

With a coupling one can typically associate two conditional distributions.
In our situation, a~coupling $M_n$
satisfying~\eqref{eq:symm_IPS_coupling}
leads to two distributions
on $I_n$~\eqref{eq:IPS_interval_for_xn}
which we denote by
\begin{gather*}
	p_n=\mathrm{Law}\big(x_n^{s_n{\boldsymbol \nu}}(t)\,|\, \mathbf{x}^{{\boldsymbol \nu}}(t)\big)
	\qquad \textnormal{and}\qquad
	p'_n=\mathrm{Law}\big(x_n^{{\boldsymbol \nu}}(t)\,|\, \mathbf{x}^{s_n{\boldsymbol \nu}}(t)\big).
\end{gather*}
Indeed, under, say, $p_n$
it suffices to specify only the
conditional distribution of $x_n^{s_n{\boldsymbol \nu}}(t)$,
as all the other locations in $x^{s_n{\boldsymbol \nu}}(t)$
stay the same.
Thus, a coupling $M_n$ satisfying~\eqref{eq:symm_IPS_coupling}
is determined by either $p_n$ or $p_n'$.

In the particular example $M_n^{\mathrm{indep}}$
\eqref{eq:symm_IPS_indep_coupling}, the distribution
$p_n$ simply corresponds to forgetting the previous location of
$x_n^{{\boldsymbol \nu}}(t)$, and selecting independently
the new particle $x_n^{s_n{\boldsymbol \nu}}(t)\in I_n$ (according
to the distribution with the parameters $s_n{\boldsymbol \nu}$)
given the remaining configuration
$x_{\hat n}^{s_n{\boldsymbol \nu}}(t)=x_{\hat n}^{{\boldsymbol \nu}}(t)$.
This distribution $p_n$ corresponding to $M_n^{\mathrm{indep}}$
can be quite complicated
as it may depend on the whole remaining configuration $\mathbf{x}_{\hat n}^{{\boldsymbol \nu}}(t)$.
This dependence may also
nontrivially incorporate the time parameter~$t$.

In this paper we
describe specific integrable parameter-symmetric particle systems
for which there exist much simpler conditional probabilities $p_n$
or $p'_n$.
Let us give a definition clarifying what we mean here by ``simpler'':
\begin{Definition}
	\label{def:local_transition}
	The conditional probability $p_n$ is said to be
	local
	if
	$p_n=\mathrm{Law}\big(x_n^{s_n{\boldsymbol \nu}}(t)\,|\, \mathbf{x}^{{\boldsymbol \nu}}(t)\big)$
	depends only on $n$, ${\boldsymbol \nu}$, and three particle locations
	$x_{n+1}^{{\boldsymbol \nu}}(t)$, $x_{n}^{{\boldsymbol \nu}}(t)$, $x_{n-1}^{{\boldsymbol \nu}}(t)$.
	The definition for $p_n'$ is analogous.

	We will interpret the local conditional probability $p_n$ as a Markov operator.
	When applied, $p_n$~leads to
	a random move
	$x_{n}^{{\boldsymbol \nu}}(t)\to x_{n}^{s_n{\boldsymbol \nu}}(t)$
	given
	$x_{n+1}^{{\boldsymbol \nu}}(t)$, $x_{n-1}^{{\boldsymbol \nu}}(t)$.
	In distribution the application of~$p_n$
	is equivalent to
	the swapping of the parameters $\nu_n \leftrightarrow \nu_{n+1}$.
	Due to this interpretation, we~will call $p_n$
	the (\emph{Markov}) \emph{swap operator}.
\end{Definition}

In the examples we consider, swap operators
will also be
independent of $t$.

\begin{Remark}
	Typically,
	only one of the probabilities $p_n$ and
	$p_n'$ can be local (and thus correspond to a swap operator). Indeed, assuming that
	$p_n$ is local, we
	can write
	\begin{gather*}
		\begin{split}
			&
			p'_n\big(x_n^{{\boldsymbol \nu}}(t)=y_n\,|\, \mathbf{x}^{s_n{\boldsymbol \nu}}(t)=\mathbf{z}\big)
			\\&\hspace{20pt}=
			p_n
			\big(x_n^{s_n{\boldsymbol \nu}}(t)=z_n\,|\,
			x_{n+1}^{{\boldsymbol \nu}}(t)=y_{n+1},x_n^{{\boldsymbol \nu}}(t)=y_n,x_{n-1}^{{\boldsymbol \nu}}(t)=y_{n-1}\big)
			\,\frac{P(\mathbf{x}^{{\boldsymbol \nu}}(t)=\mathbf{y})}{P(\mathbf{x}^{s_n{\boldsymbol \nu}}(t)=\mathbf{z})},
		\end{split}
	\end{gather*}
	where $\mathbf{y}_{\hat n}=\mathbf{z}_{\hat n}$,
	and we also assume that the probability in the denominator is nonzero.
	If one wants $p_n'$ to be local, too,
	it is necessary that the
	ratio of the probabilities
	$\frac{P(\mathbf{x}^{{\boldsymbol \nu}}(t)=\mathbf{y})}{P(\mathbf{x}^{s_n{\boldsymbol \nu}}(t)=\mathbf{z})}$
	(in which~$\mathbf{y}$,~$\mathbf{z}$ differ only by the location of the $n$-th particle)
	depends only on the four particle locations
	$x_{n+1}^{{\boldsymbol \nu}}(t)$, $x_n^{{\boldsymbol \nu}}(t)$,
	$x_n^{s_n{\boldsymbol \nu}}(t)$, $x_{n-1}^{{\boldsymbol \nu}}(t)$.
	This (quite strong)
	condition on the ratio of the probabilities does not hold
	for the particle systems considered in the present paper.
	(In particular,
	using the explicit Rakos--Sch\"utz formula
	\cite{rakos2005bethe}
	expressing
	transition probabilities in TASEP with
	particle-dependent speeds as determinants
	one can check
	that the condition fails for the usual TASEP.)
\end{Remark}

\section[Swap operators for q-Hahn TASEP]{Swap operators for $\boldsymbol q$-Hahn TASEP}\label{sec:qhahn_TASEP}

{\sloppy
In this section we examine the parameter symmetry and swap operators for the
$q$-Hahn TASEP~\cite{Povolotsky2013}.
A multiparameter version of the process
preserving its integrability is due to~\cite{BorodinPetrov2016inhom}.

}

\subsection[The q-deformed beta-binomial distribution]
{The $\boldsymbol q$-deformed beta-binomial distribution}\label{sub:phi_distribution}

We first recall the definition and properties of
the $q$-deformed beta-binomial distribution
$\varphi_{q,\mu,\nu}$ from~\cite{Corwin2014qmunu, Povolotsky2013}.
We use the standard notation for the $q$-Pochhammer symbol
$(x;q)_k=(1-x)\times(1-qx)\cdots\big(1-q^{k-1}x\big)$, $k\in \mathbb{Z}_{\ge1}$
(by agreement, $(x;q)_0=1$).

Everywhere throughout the paper we
assume that the main parameter $q$ is between $0$ and $1$.

\begin{Definition}
	\label{def:phi_distribution}
	For $m\in \mathbb{Z}_{\ge0}$, consider the following distribution on
	$\left\{ 0,1,\ldots,m \right\}$:
	\begin{gather*}
		\varphi_{q,\mu,\nu}(j\,|\, m)=
		\mu^j\,\frac{(\nu/\mu;q)_j(\mu;q)_{m-j}}{(\nu;q)_m}
		\frac{(q;q)_m}{(q;q)_j(q;q)_{m-j}},
		\qquad
		0\le j\le m.
	\end{gather*}
	When $m=+\infty$, extend the definition as
	\begin{gather*}
		\varphi_{q,\mu,\nu}(j\,|\, \infty)=
		\mu^j\frac{(\nu/\mu;q)_j}{(q;q)_j}\frac{(\mu;q)_\infty}{(\nu;q)_\infty},
		\qquad
		j\in \mathbb{Z}_{\ge0}.
	\end{gather*}
	The distribution depends on $q$
	and two other parameters~$\mu,\nu$.
\end{Definition}

When $0\le \mu\le 1$ and $\nu\le \mu$, the weights $\varphi_{q,\mu,\nu}(j\,|\, m)$
are nonnegative.\footnote{These conditions do not
exhaust the full range of $(q,\mu,\nu)$
for which the weights are nonnegative.
See, e.g., \cite[Section 6.6.1]{BorodinPetrov2016inhom}
for additional families of parameters leading to nonnegative weights.}
They also sum to~one:
\begin{gather*}
\sum_{j=0}^{m}\varphi_{q,\mu,\nu}(j\,|\, m)=1,\qquad m\in\{ 0,1,\ldots \}
\cup\{ +\infty \}.
\end{gather*}
We will need two other properties of the weights given in the next
two lemmas.

\begin{Lemma}[\cite{Barraquand_qhahn_2014, Corwin2014qmunu}]
	\label{lemma:phi_symmetry}
	The weights satisfy a symmetry property: for all
	$m,y\in \mathbb{Z}_{\ge0}$ we have
	\begin{gather*}
		\sum_{j=0}^{m}q^{jy}\varphi_{q,\mu,\nu}(j\,|\, m)
		=\sum_{k=0}^{y} q^{km}\varphi_{q,\mu,\nu}(k\,|\, y).
	\end{gather*}
	Similarly, for all $y\in \mathbb{Z}_{\ge0}$, we have
	\begin{gather*}
		\sum_{j=0}^{\infty}q^{jy}\varphi_{q,\mu,\nu}(j\,|\, \infty)=
		\varphi_{q,\mu,\nu}(0\,|\, y).
	\end{gather*}
\end{Lemma}

Define the following difference operator:
\begin{gather}
	\label{eq:nabla_mu_nu}
	(\nabla_{\mu,\nu}f)(n):=\frac{\mu-\nu}{1-\nu}\,f(n-1)+\frac{1-\mu}{1-\nu}\,f(n).
\end{gather}

The next statement is a key property of the $q$-deformed beta binomial distribution
which allows to simplify the action of certain operators defined through $\varphi_{q,\mu,\nu}$
on functions satisfying special boundary conditions. This is a manifestation of the
connection of $\varphi_{q,\mu,\nu}$ to the coordinate Bethe ansatz, as
developed in~\cite{Povolotsky2013}.
\begin{Lemma}
	\label{lemma:phi_duality}
	Fix parameters $\nu_i\in(0,1)$, $i\in \mathbb{Z}$.
	Let a function $f(n_1,\ldots,n_m )$ from
	$\mathbb{Z}^{m}$ to $\mathbb{C}$ satisfy the following
	two-body
	boundary conditions
\begin{gather}
\frac{\nu_{n_i}(1-q)}{1-q\nu_{n_i}}\,f(n_1,\ldots,n_{i}-1,n_{i+1}-1,\ldots ,n_m)
+\frac{q-\nu_{n_i}}{1-q\nu_{n_i}}\,f(n_1,\ldots,n_i,n_{i+1}-1,\ldots ,n_m)\nonumber
\\ \qquad
{}	+\frac{1-q}{1-q\nu_{n_i}}\,f(n_1,\ldots,n_i,n_{i+1},\ldots ,n_m)
-f(n_1,\ldots,n_i-1,n_{i+1},\ldots ,n_m)=0
		\label{eq:q_hahn_boundary_conditions}
		\end{gather}
	for all $\vec n\in \mathbb{Z}^m$
	such that for some $i \in \left\{ 1,\ldots,m \right\}$,
	$n_i=n_{i+1}$. $($In~\eqref{eq:q_hahn_boundary_conditions},
	only the $i$-th and the $(i+1)$-st components of $\vec n$
	are changed.$)$ Then we have
	\begin{gather}
		\label{eq:quantum_binomial_variable}
		\prod_{i=1}^{m}[\nabla_{\mu,\nu_{n}}]_i \,
		f(\underbrace{n,n,\ldots,n}_m)
		=
		\sum_{j=0}^{m}
		\varphi_{q,\mu,\nu_{n}}(j\,|\, m)\,
		f(\underbrace{n,\ldots,n }_{m-j},\underbrace{n-1,\ldots,n-1 }_{j}).
	\end{gather}
	Here $[\nabla_{\mu,\nu_{n}}]_i$ is the operator~\eqref{eq:nabla_mu_nu}
	applied in the $i$-th variable.
\end{Lemma}
\begin{proof}
	This is based on the quantum (noncommutative) binomial
	\cite[Theorem 1]{Povolotsky2013}, and is a~straightforward generalization of the
	equivalence of the free and true evolution equations
	\cite[Proposition 1.8]{Corwin2014qmunu}.
	The only difference here is that $\nu_i$'s are allowed to vary.
	However, as the application of the quantum binomial
	result depends only on the parameter $\nu_n$ associated
	to the particular $n$ in~\eqref{eq:quantum_binomial_variable},
	we see that the claim readily holds.
\end{proof}

\subsection[Multiparameter q-Hahn TASEP]
{Multiparameter $\boldsymbol q$-Hahn TASEP}\label{sub:qhahn_defn}

Here we recall the particle-inhomogeneous version of
the $q$-Hahn TASEP from \cite[Section 6.6]{BorodinPetrov2016inhom}.
Let
\begin{gather*}
	\nu_i\in(0,1),\qquad i\in \mathbb{Z}_{\ge1},\qquad
	\gamma\in\bigl[1,\sup\nolimits_{i}\nu_i^{-1}\bigr]
\end{gather*}
be parameters.
To make the system nontrivial, the $\nu_i$'s should be uniformly
bounded away from~$1$.

\begin{figure}[htpb]
	\centering
	\includegraphics[scale=1.1]{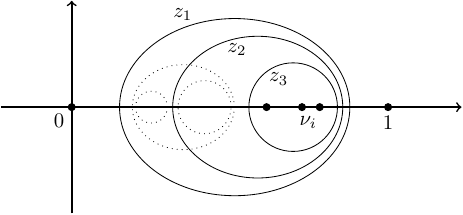}
	\caption{Possible integration contours in~\eqref{eq:qhahn_step_moments}
	for $\ell=3$. The contours for $qz_3$, $q^2z_3$, and $qz_2$ are shown dotted.}
	\label{fig:qhahn_contours}
\end{figure}

The $q$-Hahn TASEP starts from $\mathsf{step}$ and
evolves in
$\mathrm{Conf}_{\rm fin}(\mathbb{Z})$
in discrete time $t\in \mathbb{Z}_{\ge0}$.
At~each time moment, each particle $x_i$ independently
jumps to the right by $j$ with probability
\begin{gather}
	\label{eq:qhahn_one_step_jump}
	\varphi_{q,\gamma\nu_i,\nu_i}(j\,|\, x_{i-1}-x_i-1),
	\qquad
	j\in \left\{ 0,1,\ldots,x_{i-1}-x_i-1 \right\},
\end{gather}
where $x_0=+\infty$, by agreement.
See Fig.~\ref{fig:qhahn_tasep} for an illustration.
For the step initial configuration
the $q$-moments of the $q$-Hahn TASEP were
obtained in \cite[Corollary 10.4]{BorodinPetrov2016inhom}
(in the homogeneous case $\nu_i\equiv \nu$
a proof using duality and coordinate
Bethe ansatz
is due to
\cite{Corwin2014qmunu}). The $q$-moments are given in the next proposition.
\begin{Proposition}
	For any $\ell\in \mathbb{Z}_{\ge1}$
	and any $n_1\ge n_2\ge \dots\ge n_\ell\ge1$
	with the assumption that
	\begin{gather}
		\label{eq:qhahn_nu_contour_existence_assumption}
		\min_{1\le i\le n_1}
		\nu_i>q \max_{1\le i\le n_1} \nu_i
	\end{gather}
	we have for the $q$-moments of the $q$-Hahn TASEP started from $\mathsf{step}$:
	\begin{gather}
\mathbb{E}^{\mathrm{qH}(\boldsymbol\nu)}_{\mathsf{step}}
\prod_{j=1}^{\ell}q^{x_{n_j}(t)+n_j}= 
			(-1)^{\ell}q^{\frac{\ell(\ell-1)}{2}}
			\oint\frac{{\rm d}z_1}{2\pi\mathbf{i}}
			\cdots
			\oint\frac{{\rm d}z_\ell}{2\pi\mathbf{i}}
			\prod_{1\le A<B\le \ell}\frac{z_A-z_B}{z_A-qz_B}\nonumber
\\ \hphantom{\mathbb{E}^{\mathrm{qH}(\boldsymbol\nu)}_{\mathsf{step}}
\prod_{j=1}^{\ell}q^{x_{n_j}(t)+n_j}= }
{}\times\prod_{i=1}^{\ell}\left(\bigg( \frac{1-\gamma z_i}{1-z_i} \bigg)^t
\frac{1}{z_i(1-z_i)}\prod_{j=1}^{n_i}\frac{1-z_i}{1-z_i/\nu_j} \right).
\label{eq:qhahn_step_moments}
\end{gather}
	Here the integration contours are positively oriented simple
	closed curves which are
	$q$-nested around $\{\nu_j\}_{j=1,\ldots,n_1 }$
	$($that is, each contour encircles the $\nu_j$'s, and, moreover, the $z_A$ contour
	encircles each $qz_B$ contour, $B>A)$
	and leave $0$ and $1$ outside. See
	Fig.~$\ref{fig:qhahn_contours}$ for an illustration.
\end{Proposition}

\begin{Remark}
	\label{rmk:inhom_time}
	Together with particle-dependent inhomogeneity governed by
	the parameters $\nu_i$, one can make the parameter $\gamma$ time-dependent.
	That is, at each time step $t-1\to t$, the jumping distribution
	\eqref{eq:qhahn_one_step_jump} can be replaced by
	$\varphi_{q,\gamma_t\nu_i,\nu_i}(j\,|\, x_{i-1}-x_i-1)$.
	The moment formula~\eqref{eq:qhahn_step_moments} continues to hold
	when modified by replacing the term $\bigl( \frac{1-\gamma z_i}{1-z_i} \bigr)^t$
	with
	$\prod_{l=1}^t\frac{1-\gamma_l z_i}{1-z_i}$. The main result
	of this section (Theorem~\ref{thm:qhahn_swap} below)
	also holds in this generality,
	but for simplicity we will continue to assume that $\gamma$
	does not depend on $t$.
\end{Remark}

Since $0<q<1$ and we start from $\mathsf{step}$, the random variables
$\prod_{j=1}^{\ell}q^{x_{n_j}(t)+n_j}$ are all between~$0$ and~$1$.
Because the moment problem
for bounded random variables admits a unique solution,
the $q$-moments~\eqref{eq:qhahn_step_moments}
uniquely determine the joint distribution of
all the
$q$-Hahn TASEP particles $\{x_i(t)\}_{i\in \mathbb{Z}_{\ge1}}$ at each fixed time moment.
This implies the following statement:

\begin{Proposition}
	The multiparameter $q$-Hahn TASEP
	started from the step initial configuration
	is a parameter-symmetric particle system in the sense of
	Definition~$\ref{def:symm_IPS}$.
	Moreover, the distribution of each $x_n(t)$
	depends on the parameters $\nu_1,\ldots,\nu_n $ in a symmetric way.
\end{Proposition}

Denote the right-hand of~\eqref{eq:qhahn_step_moments}
by $f(n_1,\ldots,n_\ell )$, where now $n_i\in \mathbb{Z}$ are not necessarily ordered.
Notice that if $n_\ell=0$, the integrand has no poles inside the $z_\ell$ (i.e., the smallest)
integration contour.
Therefore,
$f(n_1,\ldots,n_{\ell-1},0 )=0$.
The following lemma will be employed in the next section.

\begin{Lemma}
	\label{lemma:boundary_conditions}
	The function $f(n_1,\ldots,n_\ell )$ on $\mathbb{Z}^{\ell}$ defined before the lemma
	satisfies the two-body boundary conditions
	\eqref{eq:q_hahn_boundary_conditions}.
\end{Lemma}
\begin{proof}
	This statement essentially appears in
	\cite{Corwin2014qmunu}, see also~\cite{BorodinCorwinSasamoto2012}.
	Its proof is rather short so we reproduce it here.
	When $n_i=n_{i+1}$ (denote them both by $n$),
	the part of the integrand in~\eqref{eq:qhahn_step_moments} depending on
	$z_i$, $z_{i+1}$ contains
	\begin{gather*}
		\frac{z_i-z_{i+1}}{z_i-qz_{i+1}}
		\prod_{j=1}^{n}\frac{(1-z_i)(1-z_{i+1})}{(1-z_i/\nu_j)(1-z_{i+1}/\nu_j)}.
	\end{gather*}
	The left-hand side of the boundary conditions
	\eqref{eq:q_hahn_boundary_conditions} for our function $f$
	is an integral over the contours as in Fig.~\ref{fig:qhahn_contours}, where the integrand now contains
	\begin{gather*}
	\frac{z_i-z_{i+1}}{z_i-qz_{i+1}}
	\prod_{j=1}^{n-1}\frac{(1-z_i)(1-z_{i+1})}{(1-z_i/\nu_j)(1-z_{i+1}/\nu_j)}
	\\ \hphantom{\frac{z_i-z_{i+1}}{z_i-qz_{i+1}}}
	{}\times
	\bigg(\frac{\nu_n(1\!-\!q)}{1\!-\!q \nu_n}+\frac{q\!-\!\nu_n}{1\!-\!q \nu_n}\frac{1\!-\!z_i}{1\!-\!z_i/\nu_n}
	+\frac{1\!-\!q}{1\!-\!q \nu_n}\frac{1\!-\!z_i}{1\!-\!z_i/\nu_n}\frac{1\!-\!z_{i+1}}{1\!-\!z_{i+1}/\nu_n}
	-\frac{1\!-\!z_{i+1}}{1\!-\!z_{i+1}/\nu_n}\bigg)
	\\ \hphantom{\frac{z_i-z_{i+1}}{z_i-qz_{i+1}}}
{}=\frac{\nu_n(1-\nu_n)^2}{1-q\nu_n}
	\frac{z_i-z_{i+1}}{(z_i-\nu_n)(z_{i+1}-\nu_n)}
	\prod_{j=1}^{n-1}\frac{(1-z_i)(1-z_{i+1})}{(1-z_i/\nu_j)(1-z_{i+1}/\nu_j)},
	\end{gather*}
	where the important observation is that
	the denominator $z_i-qz_{i+1}$ has
	canceled out. Now the contour for $z_{i}$ can be deformed
	(without picking any residues) to coincide with the contour for $z_{i+1}$.
	However, thanks to the factor $z_i-z_{i+1}$, the integrand is antisymmetric in
	$z_i$, $z_{i+1}$. Therefore, the whole integral vanishes, as desired.
\end{proof}

\subsection[Markov swap operators for q-Hahn TASEP]
{Markov swap operators for $\boldsymbol q$-Hahn TASEP}
\label{sub:cond_distr_qhahn}

Here we prove that the $q$-Hahn TASEP
admits a local conditional distribution
corresponding to the permutation $s_n=(n,n+1)$
when the parameters satisfy $\nu_{n+1}<\nu_n$ before their swap.
This leads to the Markov swap operator which we define now.
Fix $n\in \mathbb{Z}_{\ge1}$,
and let
\begin{gather}
	\label{eq:qhahn_transition_prob}
	p_n^{\mathrm{qH}}
	(x_n'\,|\, x_{n+1},x_n,x_{n-1}):=
	\varphi_{q,\frac{\nu_{n+1}}{\nu_n},\nu_{n+1}}(x_n'-x_{n+1}-1\,|\, x_n-x_{n+1}-1).
\end{gather}
Observe that this probability does not depend on $x_{n-1}$.
The condition $\nu_{n+1}<\nu_n$ ensures that
the swap operator $p_n^{\mathrm{qH}}$
\eqref{eq:qhahn_transition_prob}
has nonnegative probability weights.
\begin{Theorem}[Theorem~\ref{thm:intro_qhahn_swap} in Introduction]
	\label{thm:qhahn_swap}
	Let $\mathbf{x}^{\boldsymbol\nu}(t)$ be the $q$-Hahn TASEP with parameters
	$\boldsymbol\nu=\{\nu_i\}_{i\in \mathbb{Z}_{\ge1}}$,
	started from $\mathsf{step}$.
	Fix $n\in \mathbb{Z}_{\ge1}$ and assume that
	$\nu_{n+1}<\nu_n$.
	Replace $x_n(t)$ by a~random $x_n'(t)$ coming from the Markov swap operator
	$p_n^{\mathrm{qH}}$~\eqref{eq:qhahn_transition_prob}.
	Then the new configuration is distributed as the $q$-Hahn TASEP $\mathbf{x}^{s_n \boldsymbol\nu}(t)$
	with swapped parameters.
\end{Theorem}
\begin{proof}
	We will prove this theorem by applying
	$p_n^{\mathrm{qH}}$
	in the $q$-moment formula.
	Since the $q$-mo\-ments uniquely determine
	the distribution, this computation will imply the claim.
	
	Fix integers $\ell,\ell',a,b\ge0$ and define
	\begin{gather}
	\label{eq:n_a_b_notation}
	\vec n=(n_1,\ldots,n_k ):=(m_1,\ldots,m_\ell,\underbrace{n+1,\ldots,n+1 }_a,
 \underbrace{n,\ldots,n }_b,m_1',\ldots,m'_{\ell'}),
	\end{gather}
	where $m_1\ge \dots\ge m_\ell>n+1$,
	$n>m_1'\ge \dots\ge m'_{\ell'}\ge1$, and
	$k=\ell+a+b+\ell'$.
	Assume that the parameters $\nu_i$ satisfy
	the contour existence condition
	\eqref{eq:qhahn_nu_contour_existence_assumption}. (In the end
	of the proof we will drop this assumption.)
	It suffices to show that
	\begin{gather}
		\label{eq:expectation_Markov_chain_action}
		\mathbb{E}^{\mathrm{qH}(\boldsymbol\nu)}_{\mathsf{step}}
		\bigg(
			\sum_{x_n'=x_{n+1}(t)+1}^{x_n(t)}
			p_n^{\mathrm{qH}}
			(x_n'\,|\, x_{n+1}(t),x_n(t),x_{n-1}(t))
			\,q^{b(x_n'+n)}
			\prod_{\substack{j=1 \\n_j\ne n}}^{k}q^{x_{n_j}(t)+n_j}
		\biggr)
	\end{gather}
	(where the expectation
	$\mathbb{E}^{\mathrm{qH}(\boldsymbol\nu)}_{\mathsf{step}}$
	is taken with the
	parameters before the swap),
	is equal to the expectation
	\begin{gather}
		\label{eq:expectation_desired_for_qhahn}
		\mathbb{E}^{\mathrm{qH}(s_n\boldsymbol\nu)}_{\mathsf{step}}
		\prod_{j=1}^{k}q^{x_{n_j}(t)+n_j}
	\end{gather}
	with the swapped parameters.
	Indeed,
	the sum over $x_n'$
	in
	\eqref{eq:expectation_Markov_chain_action}
	corresponds to the action of the swap operator
	$p_n^{\mathrm{qH}}$
	on
	$\prod_{j=1}^{k}q^{x_{n_j}(t)+n_j}$ viewed as a function of $\{x_i(t)\}$.
	We thus need to show that the expectation of the result
	with respect to the original parameters
	leads to
	the formula with the swapped parameters.

	We now start from~\eqref{eq:expectation_Markov_chain_action}, and
	in the rest of the proof omit the dependence on $t$ for shorter notation.
	First, we use the symmetry property (Lemma~\ref{lemma:phi_symmetry})
	to write for the part of the sum in~\eqref{eq:expectation_Markov_chain_action}
	involving $x_n$, $x_{n+1}$:
	\begin{gather}
	\sum_{x_n'=x_{n+1}+1}^{x_n}
	\varphi_{q,\frac{\nu_{n+1}}{\nu_n},\nu_{n+1}}(x_n'-x_{n+1}-1\,|\, x_n-x_{n+1}-1)
	\,q^{a(x_{n+1}+n+1)+b(x_n'+n)}\nonumber
	\\[-.5ex] \qquad
{}=	\sum_{x_n'=x_{n+1}+1}^{x_n}	\!\!\! q^{b(x_n'-x_{n+1}-1)}
	\varphi_{q,\frac{\nu_{n+1}}{\nu_n},\nu_{n+1}}
	(x_n'-x_{n+1}-1\,|\, x_n-x_{n+1}-1)	\,q^{(a+b)(x_{n+1}+n+1)}\nonumber
	\\[-.5ex] \qquad
{}=	\sum_{r=0}^{b}	q^{r(x_n-x_{n+1}-1)}	\varphi_{q,\frac{\nu_{n+1}}{\nu_n},\nu_{n+1}}
	(r\,|\, b)	\,q^{(a+b)(x_{n+1}+n+1)}\nonumber
	\\[-.5ex] \qquad
{}=	\sum_{r=0}^{b}	\varphi_{q,\frac{\nu_{n+1}}{\nu_n},\nu_{n+1}}
	(r\,|\, b)	\,q^{(a+b-r)(x_{n+1}+n+1)+r(x_n+n)}.
	\label{eq:duality_computation_for_ptransition}
	\end{gather}

	We thus need to compute
	\begin{gather}
		\label{eq:thm_qhahn_swap_proof1}
		\sum_{r=0}^{b}
		\varphi_{q,\frac{\nu_{n+1}}{\nu_n},\nu_{n+1}}
		(r\,|\, b)
		\,\mathbb{E}^{\mathrm{qH}(\boldsymbol\nu)}_{\mathsf{step}}
		\biggl(\prod_{j=1}^{k}
		q^{x_{n_j(r)}+n_j(r)}\biggr),
	\end{gather}
	where the vector
	$\vec n(r)=(n_1(r),\ldots,n_k(r) )$
	is as in~\eqref{eq:n_a_b_notation},
	but with $(a,b)$ replaced by $(a+b-r,r)$.
	The expectation in~\eqref{eq:thm_qhahn_swap_proof1}
	is given by the contour integral as in the right-hand side of
	\eqref{eq:qhahn_step_moments}.
	Recall the notation $f(\vec n)$ for this integral, where now $\vec n\in \mathbb{Z}^k$,
	and the components of $\vec n$ are not necessarily ordered.
	By Lemma~\ref{lemma:boundary_conditions}, this function $f$ satisfies the two-body boundary conditions.
	Thus,~\eqref{eq:thm_qhahn_swap_proof1} can be rewritten by Lemma~\ref{lemma:phi_duality}
	as (recall notation~\eqref{eq:nabla_mu_nu} for
	the operator
	$\nabla_{\mu,\nu}$):
	\begin{gather*}
		\prod_{j=1}^{b}\left[ \nabla_{\frac{\nu_{n+1}}{\nu_n},\nu_{n+1}} \right]_{\ell+a+j}\,
		f(m_1,\ldots,m_\ell,\underbrace{n+1,\ldots,n+1 }_{a+b},m_1',\ldots,m'_{\ell'}).
	\end{gather*}
	Observe that now
	each of the difference operators
	$[ \nabla_{\frac{\nu_{n+1}}{\nu_n},\nu_{n+1}} ]_{\ell+a+j}$
	can be applied independently
	inside the integral. We thus have for every variable $w=z_{\ell+a+j}$,
	$j=1,\ldots,b $:
	\begin{gather*}
	[\nabla_{\frac{\nu_{n+1}}{\nu_n},\nu_{n+1}}]_{\ell+a+j}
	\prod_{i=1}^{n+1}\frac{1\!-\!w}{1\!-\!w/\nu_i}
	=	\bigg( \frac{\nu_{n+1}/\nu_n\!-\!\nu_{n+1}}{1\!-\!\nu_{n+1}}+
	\frac{1\!-\!\nu_{n+1}/\nu_n}{1\!-\!\nu_{n+1}}\frac{1\!-\!w}{1\!-\!w/\nu_{n+1}} \bigg)
	\prod_{i=1}^{n}\frac{1\!-\!w}{1\!-\!w/\nu_i}
	\\ \hphantom{[\nabla_{\frac{\nu_{n+1}}{\nu_n},\nu_{n+1}}]_{\ell+a+j}
	\prod_{i=1}^{n+1}\frac{1\!-\!w}{1\!-\!w/\nu_i}}
{}=	\frac{1-w/\nu_n}{1-w/\nu_{n+1}}	\prod_{i=1}^{n}\frac{1-w}{1-w/\nu_i}
	=\prod_{\substack{i=1\\i\ne n}}^{n+1}\frac{1-w}{1-w/\nu_i}.
	\end{gather*}
	We see that the
	resulting integral coming from~\eqref{eq:expectation_Markov_chain_action}
	contains, for each variable
	$z_{\ell+a+j}$ corresponding to $n_{\ell+a+j}=n$ in~\eqref{eq:n_a_b_notation}, the product over the
	parameters $(\nu_1,\ldots,\nu_{n-1},\nu_{n+1} )=s_n(\nu_1,\ldots,\nu_{n} )$.
	Therefore, this integral is equal to the
	expectation~\eqref{eq:expectation_desired_for_qhahn}
	with the swapped parameters $s_n\boldsymbol\nu$, as desired.

	It remains to show that we can drop the contour
	existence assumption~\eqref{eq:qhahn_nu_contour_existence_assumption}.
	The preceding argument implies that
	under~\eqref{eq:qhahn_nu_contour_existence_assumption} (with fixed $x_{n+1},x_n',x_{n-1}$),
	\begin{gather}
		\sum_{x_n=x_n'}^{x_{n-1}-1}
		P_{\boldsymbol\nu}(\ldots,x_{n+1},x_n,x_{n-1},\ldots )\,
		\varphi_{q,\frac{\nu_{n+1}}{\nu_n},\nu_{n+1}}(x_n'-x_{n+1}-1\,|\, x_n-x_{n+1}-1)\nonumber
		\\ \qquad
{}=		P_{s_n\boldsymbol\nu}(\ldots, x_{n+1},x_{n}',x_{n-1},\ldots ),
		\label{eq:analytic_continuation}
	\end{gather}
	where $P_{\boldsymbol\nu}$, $P_{s_n\boldsymbol\nu}$ denote the
	$q$-Hahn probability distributions
	with the corresponding parameters
	at some fixed time $t\in \mathbb{Z}_{\ge0}$.
	
	If $n\ge2$, the sum in the left-hand side of
	\eqref{eq:analytic_continuation}
	is finite, and each probability $P_{\boldsymbol\nu}$,
	$P_{s_n\boldsymbol\nu}$ is a rational function of $\nu_2,\nu_3,\ldots $
	(note that since $\ldots,x_{n+1},x_n',x_{n-1},\ldots, $ are
	fixed, only finitely many of the $\nu_i$'s
	enter~\eqref{eq:analytic_continuation}).
	The dependence on $\nu_1$ is also rational after canceling out the common
	factor $\frac{(\gamma \nu_1;q)_\infty^t}{(\nu_1;q)_\infty^t}$
	from both sides. Therefore, identity~\eqref{eq:analytic_continuation}
	between rational functions in $\nu_i$ can be analytically continued,
	and the assumption
	\eqref{eq:qhahn_nu_contour_existence_assumption} can be dropped.

	For $n=1$, the sum in the left-hand side of~\eqref{eq:analytic_continuation}
	becomes infinite. Remove the
	common factor
	$\frac{(\gamma \nu_1;q)_\infty^t}{(\nu_1;q)_\infty^t}$
	from both sides again,
	then the coefficients by each power
	$\gamma^m$, $m\in \mathbb{Z}_{\ge0}$,
	become rational functions in $\nu_i$, $i=1,2,\ldots $. Therefore,
	we can again analytically continue identity~\eqref{eq:analytic_continuation}
	and drop the assumption~\eqref{eq:qhahn_nu_contour_existence_assumption}.
	This completes the proof.
\end{proof}

\begin{Remark}
	When $\nu_n=\nu_{n+1}$,
	we have from~\eqref{eq:qhahn_transition_prob}
	that $p_n^{\mathrm{qH}}(x_n'\,|\, x_{n+1},x_n,x_{n-1})=\mathbf{1}_{x_n'=x_n}$
	(where $\mathbf{1}_{\cdots}$ stands for the indicator).
	Therefore, the swap operator
	reduces to the identity map, which is appropriate since
	for $\nu_n=\nu_{n+1}$
	there is nothing to swap.

	If $\nu_n<\nu_{n+1}$,
	formula~\eqref{eq:qhahn_transition_prob} for $p_n^{\mathrm{qH}}$ also
	makes sense, but some of these probability weights become negative.
	One can check that
	all algebraic manipulations in the proof of Theorem~\ref{thm:qhahn_swap}
	are still valid for $\nu_n<\nu_{n+1}$,
	but now they do not correspond to actual stochastic objects.
	This is the reason for the restriction $\nu_n>\nu_{n+1}$ in
	Theorem~\ref{thm:qhahn_swap}.
\end{Remark}

\subsection[Duality for the q-Hahn swap operator]
{Duality for the $ \boldsymbol q$-Hahn swap operator}\label{sub:duality}

Here let us recall the Markov duality relation for the
$q$-Hahn TASEP from~\cite{Corwin2014qmunu}.
We will heavily use duality in Section~\ref{sec:stat_qtasep} below.

\begin{figure}[htpb]
	\centering
	\includegraphics{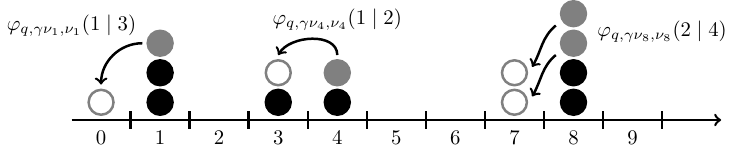}
	\caption{Configuration of particles
	$\vec n=(8,8,8,8,4,4,3,1,1,1)\in \mathbb{W}^{10}$,
	and a possible one-step transition
	in the $q$-Hahn Boson process. The particles of
	$\vec n$ which jump are solid gray, and
	their new locations are not filled.}
	\label{fig:Boson_space}
\end{figure}

Fix $\ell\ge1$ and
let
\begin{gather}
	\label{eq:W_space}
	\mathbb{W}^{\ell}:=\{\vec n=(n_1\ge \dots\ge n_\ell\ge0 ),\
	n_i\in \mathbb{Z}\}.
\end{gather}
We interpret elements of $\mathbb{W}^{\ell}$ as
$\ell$-particle configurations in $\mathbb{Z}_{\ge0}$,
where multiple particles per site are allowed. Namely,
for each $i=1,\ldots,\ell $, we put one particle at
the location $n_i\in \mathbb{Z}_{\ge0}$.
See Fig.~\ref{fig:Boson_space} for an illustration.

Define the \emph{duality functional}
on the product space
$\mathrm{Conf}_{\rm fin}(\mathbb{Z})\times\mathbb{W}^{\ell}$
as follows:
\begin{gather}
	\label{eq:duality_functional}
	H(\mathbf{x},\vec n):=\begin{cases}
		\prod\limits_{i=1}^{\ell}q^{x_{n_i}+n_i},& n_\ell\ge1,\\
		0,&n_\ell=0.
	\end{cases}
\end{gather}

Let $\mathcal{T}^{\mathrm{qH(\boldsymbol\nu)}}(\mathbf{x},\mathbf{y})$,
$\mathbf{x},\mathbf{y}\in
\mathrm{Conf}_{\rm fin}(\mathbb{Z})$,
denote the one-step Markov transition operator
of the $q$-Hahn TASEP with parameters $\boldsymbol\nu=\{\nu_i\}$ and $\gamma$.
We do not include the latter in the notation
and assume that it is fixed throughout this section.

Let
$\breve{\mathcal{T}}^{\mathrm{qH(\boldsymbol\nu)}}
(\vec n,\vec m)$, $\vec n,\vec m\in \mathbb{W}^{\ell}$,
be the one-step transition operator of a discrete time
Markov chain on $\mathbb{W}^{\ell}$
which at each time step evolves as follows.
Independently at every site $k\in \mathbb{Z}_{\ge1}$
containing, say, $y_k$ particles,
randomly select $j$ particles with probability
$\varphi_{q,\gamma\nu_k,\nu_k}(j\,|\, y_k)$,
and move them to the site $k-1$.
This Markov chain is called the (stochastic)
\emph{$q$-Hahn Boson process}.
See Fig.~\ref{fig:Boson_space} for an illustration.

\begin{Proposition}[\cite{Corwin2014qmunu}]
	With the above definitions, we have
	\begin{gather*}
		\mathcal{T}^{\mathrm{qH(\boldsymbol\nu)}}
		H(\mathbf{x},\vec n)=
		\breve{\mathcal{T}}^{\mathrm{qH(\boldsymbol\nu)}}
		H(\mathbf{x},\vec n),\qquad
		\mathbf{x}\in \mathrm{Conf}_{\rm fin}(\mathbb{Z}),
		\qquad \vec n\in \mathbb{W}^{\ell}.
	\end{gather*}
	Here the operators
	$\mathcal{T}^{\mathrm{qH(\boldsymbol\nu)}}$,
	$\breve{\mathcal{T}}^{\mathrm{qH(\boldsymbol\nu)}}$ act
	in the $\mathbf{x}$ and the $\vec n$ variables, respectively.
	Equivalently in~terms of expectations, we have
	for all
	$\mathbf{x}^0\in \mathrm{Conf}_{\rm fin}(\mathbb{Z})$,
	$\vec n^0\in \mathbb{W}^{\ell}$,
	and all times
	$t\in \mathbb{Z}_{\ge0}$:
	\begin{gather}
		\mathbb{E}^{\mathrm{qH}(\boldsymbol\nu)}_{\mathbf{x}(0)=\mathbf{x}^0}
		H\big(\mathbf{x}(t),\vec n^0\big)=
		\mathbb{E}^{\mathrm{qHBoson}(\boldsymbol\nu)}_{\vec n(0)=\vec n^0}
		H\big(\mathbf{x}^0,\vec n(t)\big).
	\end{gather}
	Here in the left-hand side the expectation is taken
	with respect to the $q$-Hahn TASEP's evolution
	starting from $\mathbf{x}^0$,
	and in the right-hand side the expectation
	is with respect to the $q$-Hahn Boson process started from
	$\vec n^0$.
\end{Proposition}

\begin{figure}[htpb]
	\centering
	\includegraphics{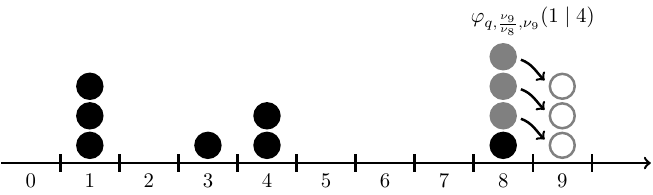}
	\caption{A possible transition under $\breve{p}_8^{\mathrm{qH}}$.}
	\label{fig:Boson_space_transition}
\end{figure}

Consider now the Markov swap operator
$p_k^{\mathrm{qH}}$~\eqref{eq:qhahn_transition_prob}
on $\mathrm{Conf}_{\rm fin}(\mathbb{Z})$,
where $k\in \mathbb{Z}_{\ge1}$ is fixed, and $\nu_k>\nu_{k+1}$.
It turns out that this operator admits a dual Markov operator
on the space~$\mathbb{W}^{\ell}$,
by means of the same duality functional $H$~\eqref{eq:duality_functional}.
Namely,
define $\breve{p}_k^{\mathrm{qH}}$
as the Markov operator which acts only on the
$k$-th location in the $q$-Boson configuration.
If the $k$-th location contains~$y_k$ particles,
then $\breve{p}_k^{\mathrm{qH}}$
randomly sends $y_k-j$ particles from location $k$
to location $k+1$, with probability
$\varphi_{q,\frac{\nu_{k+1}}{\nu_k},\nu_{k+1}}(j\,|\, y_k)$.
See Fig.~\ref{fig:Boson_space_transition} for an illustration.

\begin{Proposition}
	\label{prop:qhahn_transition_duality}
	If $\nu_k>\nu_{k+1}$,
	then we have
	\begin{gather*}
		p_k^{\mathrm{qH}}
		H(\mathbf{x},\vec n)=
		\breve{p}_k^{\mathrm{qH}}
		H(\mathbf{x},\vec n),
	\end{gather*}
	where the operator in the left-hand side acts
	on $\mathbf{x}$,
	and in the right-hand side~-- on $\vec n$.
\end{Proposition}
\begin{proof}
	This duality relation immediately follows from
	computation~\eqref{eq:duality_computation_for_ptransition}
	performed (with the help of
	Lemma~\ref{lemma:phi_symmetry})
	in the proof of Theorem~\ref{thm:qhahn_swap}.
\end{proof}

\section{Continuous time limit of repeated swaps}
\label{sec:qhahn_cont_time}

Here we consider two continuous time limits,
one of the original $q$-Hahn TASEP, and another
one of the transition probabilities $p_n^{\mathrm{qH}}$
leading to the new \emph{backward $q$-Hahn process}.
These two continuous time processes
act on the two-parameter family
$\{\mathscr{M}^{\mathrm{qH}}_{q,\nu;\mathsf{t}}\}
_{\mathsf{t}\in \mathbb{R}_{\ge0},\; 0\le \nu<1}$
of
distributions of the continuous time
$q$-Hahn TASEP (started from $\mathsf{step}$) by
continuously changing the para\-meters.

\subsection[Two expansions of the distribution phi]
{Two expansions of the distribution $\boldsymbol{\varphi_{q,\mu,\nu}}$}

Let us write down two Taylor expansions of the
$q$-deformed beta-binomial distribution
from Section~\ref{sub:phi_distribution}. Their proofs are straightforward.

\begin{Lemma}
	For $\nu\in[0,1)$
	and $m \in\{1,2,\ldots, \}\cup\{ +\infty \}$
	we have as $\varepsilon\searrow0{:}$
	\begin{gather*}
		\varphi_{q,\nu+\varepsilon,\nu}(j\,|\, m)=
		\begin{cases}
			1+O(\varepsilon),&j=0,
			\\[4pt]
			\dfrac{\nu^{j-1}}{1-q^j}
			\dfrac{(q;q)_m(\nu;q)_{m-j}}{(q;q)_{m-j}(\nu;q)_m}
			\,\varepsilon
			+O\big(\varepsilon^2\big),&1\le j\le m.
		\end{cases}
	\end{gather*}
\end{Lemma}
\begin{Lemma}
	\label{lemma:phi_expansion_for_bwd}
	For $\nu\in[0,1)$ and $m \in\{ 0,1,2,\ldots \}$
	we have as $\varepsilon\searrow0{:}$
	\begin{gather*}
		\varphi_{q,1-\varepsilon,\nu(1-\varepsilon)}(j\,|\, m)
		=
		\begin{cases}
			\dfrac{1}{1-q^{m-j}}
			\dfrac{(q;q)_m(\nu;q)_j}
			{(q;q)_j(\nu;q)_{m}}
			\,\varepsilon+O\big(\varepsilon^2\big),&0\le j\le m-1,
			\\[8pt]
			1+O(\varepsilon),&j=m.
		\end{cases}
	\end{gather*}
\end{Lemma}

Denote
\begin{gather}
	\label{eq:psi_things_notation}
	\psi_{q,\nu}(j\,|\, m):=
	\dfrac{\nu^{j-1}}{1-q^j}
	\dfrac{(q;q)_m(\nu;q)_{m-j}}{(q;q)_{m-j}(\nu;q)_m},
	\qquad
	\psi^{\bullet}_{q,\nu}(j'\,|\, m):=
	\dfrac{1}{1-q^{m-j'}}
	\dfrac{(q;q)_m(\nu;q)_{j'}}
	{(q;q)_{j'}(\nu;q)_{m}},
\end{gather}
where $1\le j\le m$ and $0\le j'\le m-1$.
Clearly, $\psi_{q,\nu}(j\,|\, m)=\nu^{j-1}\psi_{q,\nu}^{\bullet}(m-j\,|\, m)$
for all $1\le j\le m$, but it is convenient to keep these notations separate.

\subsection[Continuous time q-Hahn TASEP]
{Continuous time $\boldsymbol q$-Hahn TASEP}\label{sub:cont_time_qhahn}

Fix $q\in(0,1)$, $\nu\in[0,1)$
and consider the Poisson-type
limit of the $q$-Hahn TASEP with homogeneous parameters $\nu_i\equiv \nu$ as
\begin{gather*}
	\gamma=1+\varepsilon/\nu,	\qquad
	t=\lfloor \mathsf{t}/\varepsilon \rfloor,\qquad
	\varepsilon\searrow0,
\end{gather*}
where $t\in \mathbb{Z}_{\ge0}$ is the discrete time before the limit,
and $\mathsf{t}\in \mathbb{R}_{\ge0}$
is the scaled continuous time after the limit.
The resulting process is the \emph{continuous time $q$-Hahn TASEP}
which evolves as follows.
Starting from $\mathsf{step}$,
in continuous time $\mathsf{t}\in \mathbb{R}_{\ge0}$,
each particle $x_n(\mathsf{t})$, $n\in \mathbb{Z}_{\ge1}$, independently
jumps to the right by
$j\in \left\{ 1,2,\ldots,x_{n-1}(\mathsf{t})-x_n(\mathsf{t})-1 \right\}$
at rate
\begin{gather*}
	\psi_{q,\nu}(j\,|\, x_{n-1}(\mathsf{t})-x_{n}(\mathsf{t})-1).
\end{gather*}
(In continuous time there is at most one jump at every
instance of time.)
Here $x_0\equiv+\infty$, by agreement.
This continuous time process was considered in
\cite{barraquand2015q, Takeyama2014}.

\begin{Remark}[multiparameter continuous time $q$-Hahn TASEP]	\label{rmk:inhom_cont_time_qhahn}
	One can also consider a version of the continuous time $q$-Hahn TASEP
	in which each particle $x_n$ jumps at
	rate $\nu_n\psi_{q,\nu_n}(\cdot\,|\, \cdot)$.
	This~multiparameter deformation preserves integrability.
	To get from the multiparameter
	process to the homogeneous one described above one has to set $\nu_n\equiv \nu$
	and rescale the continuous time by the factor of $\nu$.
	For simpler notation,
	we will mostly consider the homogeneous continuous time $q$-Hahn TASEP.
	Its multiparameter generalization is needed only in the proof of~Theorem~\ref{thm:action_on_qhahn_distributions} below.
\end{Remark}

The continuous time $q$-Hahn TASEP possesses the
same $q$-moment formulas as~\eqref{eq:qhahn_step_moments},
with $\nu_i\equiv \nu$, and the replacement
\begin{gather*}
	\bigg( \frac{1-\gamma z_i}{1-z_i} \bigg)^{t}\to
	\exp\bigg\{{-}\frac{\mathsf{t}}{\nu}\frac{z_i}{1-z_i} \bigg\},
	\qquad i=1,\ldots,\ell,
\end{gather*}
inside the contour integral.

\subsection[Backward q-Hahn process]
{Backward $\boldsymbol q$-Hahn process}\label{sub:cont_time_backwards_process}

Here we define a continuous time limit of a
certain combination of the
swap operators~$p_n^{\mathrm{qH}}$~\eqref{eq:qhahn_transition_prob}.
Let us first explain the main idea.
Assume that $\nu_1>\nu_2>\nu_3>\cdots $.
By Theorem~\ref{thm:qhahn_swap},
the application of the Markov operators
$p_1^{\mathrm{qH}}$, $p_2^{\mathrm{qH}}$, \ldots
(in this order)
to a random configuration
coming from the discrete time $q$-Hahn TASEP
with parameters $\boldsymbol\nu$
is equivalent in distribution
to
the permutation
$\boldsymbol\nu\mapsto \cdots s_3s_2s_1 \boldsymbol\nu$ which
exchanges $\nu_1$ with $\nu_2$, then $\nu_1$ with $\nu_3$, and so on
(so that $\nu_1$ gets pushed all the way to infinity and essentially disappears).
Because the particle
configuration to which we apply the $p_i^{\mathrm{qH}}$'s
is densely packed to the left,
this application of the
infinite product of
the swap operators
$p_i^{\mathrm{qH}}$ is well-defined and is a one-step
Markov transition operator on $\mathrm{Conf}_{\rm fin}(\mathbb{Z})$.
Its continuous time limit as $\nu_i\to \nu$ for all $i$ will be our new backward $q$-Hahn TASEP.

Let us now make this idea precise and take the particular parameters
$\nu_i=\nu r^{i-1}$, $i\in \mathbb{Z}_{\ge0}$,
where
$\nu\in[0,1)$ and $r\in(0,1)$.
Denote by $p_n^{\mathrm{qH};(\alpha,\beta)}$, $0<\beta<\alpha$, the Markov swap operator
acting on the $n$-th particle as follows
\begin{gather*}
	p_n^{\mathrm{qH};(\alpha;\beta)}
	(x_n'\,|\, x_{n+1},x_n,x_{n-1})=
	\varphi_{q,\frac{\beta}{\alpha},\beta}(x_n'-x_{n+1}-1\,|\, x_n-x_{n+1}-1).
\end{gather*}
Then define the infinite product
\begin{gather*}
	\mathcal{B}^{\mathrm{qH}(r)}_{q,\nu}:=
	p_1^{\mathrm{qH};(\nu_2,\nu_1)}
	p_2^{\mathrm{qH};(\nu_3,\nu_1)}
	p_3^{\mathrm{qH};(\nu_4,\nu_1)}
	\cdots
	=	
	p_1^{\mathrm{qH};(\nu r,\nu)}
	p_2^{\mathrm{qH};(\nu r^2,\nu)}
	p_3^{\mathrm{qH};(\nu r^3,\nu)}
	\cdots
\end{gather*}
(these are Markov operators so their product is written as if it's the action on
measures: we first apply $p_1^{\mathrm{qH}}$, then $p_2^{\mathrm{qH}}$, and so on).
In words, under the Markov operator
$\mathcal{B}^{\mathrm{qH}(r)}_{q,\nu}$,
each particle $x_n$
jumps to the left
into a new location $x_n'\in\{x_{n+1}+1,x_{n+1}+2,\ldots,x_n \}$ chosen randomly
from the distribution
\begin{gather}
	\label{eq:phi_geometric_rates}
	\varphi_{q,r^n,\nu r^n}(x_n'-x_{n+1}-1\,|\, x_n-x_{n+1}-1).
\end{gather}
The update is sequential for $n=1,2,\ldots $, so
the new location $x_n'$ of each $x_n$ depends only on
the two old locations $x_{n+1}$, $x_n$.

\begin{Proposition}
	\label{prop:qhahn_discrete_repeated_action}
	For $\nu_i=\nu r^{i-1}$ and any $m,k\in \mathbb{Z}_{\ge0}$ we have
	\begin{gather*}
		\delta_{\mathsf{step}}\,
		\big(\mathcal{T}^{\mathrm{qH}(\boldsymbol\nu)}\big)^m
		\mathcal{B}^{\mathrm{qH}(r)}_{q,\nu}
		\mathcal{B}^{\mathrm{qH}(r)}_{q,r\nu}
		\mathcal{B}^{\mathrm{qH}(r)}_{q,r^2\nu}
		\cdots
		\mathcal{B}^{\mathrm{qH}(r)}_{q,r^{k-1}\nu}
		=
		\delta_{\mathsf{step}}\,
		\big(\mathcal{T}^{\mathrm{qH}(r^k\boldsymbol\nu)}\big)^m,
	\end{gather*}
	where
	$\mathcal{T}^{\mathrm{qH}(\boldsymbol\nu)}$ is the one-step Markov transition operator
	of the discrete time $q$-Hahn TASEP,
	$\delta_{\mathsf{step}}$ is the delta measure at the step configuration,
	and
	$r^k \boldsymbol\nu=\{\nu r^{k+i-1}\}_{i\in \mathbb{Z}_{\ge1}}$
	is the parameter sequence shifted by $k$.
\end{Proposition}
\begin{proof}
	Follows from
	Theorem~\ref{thm:qhahn_swap}.
\end{proof}

Note that
$\mathcal{T}^{\mathrm{qH}(r^k\boldsymbol\nu)}$, the $q$-Hahn TASEP
transition operator with the shifted parameter sequ\-ence $r^k\boldsymbol\nu$ from
Proposition~\ref{prop:qhahn_discrete_repeated_action},
is the same as
$\mathcal{T}^{\mathrm{qH}(\boldsymbol\nu)}\big\vert_{\nu\mapsto r^k \nu}$,
i.e., the original $q$-Hahn TASEP operator
in which $\nu$ is replaced by $r^k \nu$.

\begin{Proposition}
	\label{prop:limit_bwd_discrete_to_cont}
	In the scaling regime
	\begin{gather*}
		r=1-\varepsilon, \qquad
		k=\lfloor \tau/\varepsilon \rfloor,
		\qquad
		k'=\lfloor \tau'/\varepsilon \rfloor,
		\qquad
		\varepsilon\searrow0,
	\end{gather*}
	where $0\le \tau\le \tau'$ are scaled times,
	and
	the parameters $q,\nu$ are assumed fixed,
	we have
	\begin{gather}
		\label{eq:limit_backwards_to_cont}
		\lim_{\varepsilon\searrow0}
		\bigl(
			\mathcal{B}^{\mathrm{qH}(r)}_{q,r^k\nu}
			\mathcal{B}^{\mathrm{qH}(r)}_{q,r^{k+1}\nu}
			\mathcal{B}^{\mathrm{qH}(r)}_{q,r^{k+2}\nu}
		\cdots
		\mathcal{B}^{\mathrm{qH}(r)}_{q,r^{k'-2}\nu}
		\mathcal{B}^{\mathrm{qH}(r)}_{q,r^{k'-1}\nu}
		\bigr)
		=\mathcal{B}^{\mathrm{qH}}_{q,\nu}(\tau,\tau')
	\end{gather}
	as Markov operators acting on the space $\mathrm{Conf}_{\rm fin}(\mathbb{Z})$.
	The convergence is in the sense of the operators' matrix elements $($i.e., the strong operator
	topology$)$.
	The operators $\mathcal{B}^{\mathrm{qH}}_{q,\nu}(\tau,\tau')$ form a~continuous time,
	time-inhomogeneous
	Markov semigroup.
\end{Proposition}
The fact that the resulting continuous time Markov chain is time-inhomogeneous
will become clearer later in Section~\ref{sub:cont_time_action}
when we consider its action on the $q$-Hahn TASEP distributions.
\begin{proof}[Proof of Proposition~\ref{prop:limit_bwd_discrete_to_cont}]
	Observe that the space of left-packed configurations
	$\mathrm{Conf}_{\rm fin}(\mathbb{Z})$
	is countable, and under the Markov operators
	in both sides of
	\eqref{eq:limit_backwards_to_cont}
	the particles jump only to the left.
	Therefore, the desired limit as $\varepsilon\searrow0$
	reduces to the limit
	of finite-size Markov transition matrices.
	For the latter
	the Poisson-type limit is taken in a straightforward way.
\end{proof}

Let us now describe the
time-inhomogeneous Markov dynamics
$\mathcal{B}^{\mathrm{qH}}_{q,\nu}(\tau,\tau')$
in terms of gene\-ra\-tors.
Taking the limit $\varepsilon\searrow0$ in
the probabilities~\eqref{eq:phi_geometric_rates}
and using Lemma~\ref{lemma:phi_expansion_for_bwd}
leads to the jump rates
\begin{gather}
	\label{eq:new_new_psi_dot_jump_rates}
	n\cdot\psi^{\bullet}_{q,\nu}(x_n'-x_{n+1}-1\,|\, x_n-x_{n+1}-1)
\end{gather}
with which
each particle $x_n$
jumps to the left into $x_n'\in\{x_{n+1}+1,x_{n+1}+2,\ldots,x_n-1 \}$.
The fac\-tor~$n$ appears from the expansion $r^n=(1-\varepsilon)^n=
1-n\varepsilon+O\big(\varepsilon^2\big)$.
Denote the Markov generator with the jump rates~\eqref{eq:new_new_psi_dot_jump_rates}
by
$\mathsf{B}^{\mathrm{qH}}_{q,\nu}$.
The action of this generator is well-defined
because the configurations from $\mathrm{Conf}_{\rm fin}(\mathbb{Z})$
are densely packed to the left, so only finitely many particles
can jump in finite time.
Propositions~\ref{prop:qhahn_discrete_repeated_action} and~\ref{prop:limit_bwd_discrete_to_cont} then imply that
the semigroup and the generator are related as
\begin{gather*}
	\mathcal{B}^{\mathrm{qH}}_{q,\nu}(\tau,\tau')=
	\exp
	\bigg\{ \int_\tau^{\tau'}\mathsf{B}^{\mathrm{qH}}_{q,\nu {\rm e}^{-s}}\,{\rm d}s \bigg\},
	\qquad \qquad
	\mathsf{B}^{\mathrm{qH}}_{q,\nu}=\frac{\rm d}{{\rm d}\tau'}\Big\vert_{\tau'=0}
	\mathcal{B}^{\mathrm{qH}}_{q,\nu}(0,\tau').
\end{gather*}
In words, $\mathcal{B}^{\mathrm{qH}}_{q,\nu}(\tau,\tau')$
is the Markov transition operator from time $\tau$ to time $\tau'$
of a process under which at each time $s$
the jumps are governed by
the infinitesimal generator $\mathsf{B}^{\mathrm{qH}}_{q,\nu {\rm e}^{-s}}$.
We call the process corresponding to
$\mathcal{B}^{\mathrm{qH}}_{q,\nu}(\tau,\tau')$
the \emph{backward $q$-Hahn process}.

\begin{Remark}
	A time- and space-homogeneous version of the
	backward $q$-Hahn process
	was considered in~\cite{barraquand2015q}.
	Indeed, the
	left jumps
	with rates $\phi^L_{q,\nu}(x_n-x_n'\,|\, x_n-x_{n+1}-1)$
	in the
	$q$-Hahn asymmetric exclusion process
	from~\cite{barraquand2015q}
	coincide (up to the constant factor $L$)
	with the jumps at~rates
	$\psi^{\bullet}_{q,\nu}(x_n'-x_{n+1}-1\,|\, x_n-x_{n+1}-1)$.
	However, the spatial inhomogeneity
	of the backward $q$-Hahn process does not allow to immediately apply the
	contour integral $q$-moment formulas from
	\cite{barraquand2015q} in our situation.
\end{Remark}

\subsection{Action on distributions}
\label{sub:cont_time_action}

\begin{figure}[htpb]
	\centering
	\includegraphics{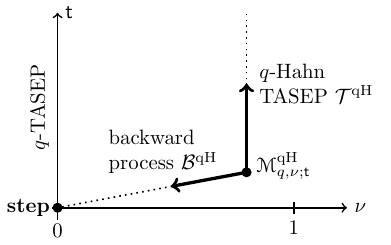}
	\caption{Action of the continuous time
	$q$-Hahn TASEP
	and the
	backward $q$-Hahn process
	on the measures
	$\mathscr{M}^{\mathrm{qH}}_{q,\nu;\mathsf{t}}$
	viewed as a two-parameter family depending
	on $\mathsf{t}$ and $\nu$.
	The vertical line $\nu=0$ corresponds to
	distributions of the $q$-TASEP,
	and on them we obtain a stationary dynamics discussed in~Section~\ref{sec:stat_qtasep} below.}
	\label{fig:two_times_qhahn_diagram}
\end{figure}

For $\mathsf{t}\in\mathbb{R}_{\ge0}$ and $\nu\in [0,1)$
denote by
$\mathscr{M}^{\mathrm{qH}}_{q,\nu;\mathsf{t}}$
the time $\mathsf{t}$
distribution of the continuous
$q$-Hahn TASEP started from $\mathsf{step}$.
The combined results of Sections~\ref{sub:cont_time_qhahn} and \ref{sub:cont_time_backwards_process}
imply the following action of the backward $q$-Hahn
process on these distributions:

\begin{Theorem}
	\label{thm:action_on_qhahn_distributions}
	We have
	for all $\nu\in[0,1)$ and $\mathsf{t},\tau\in \mathbb{R}_{\ge0}$:
	\begin{gather*}
		\mathscr{M}^{\mathrm{qH}}_{q,\nu;\mathsf{t}}\,
		\mathcal{B}^{\mathrm{qH}}_{q,\nu}(0,\tau)=
		\mathscr{M}^{\mathrm{qH}}_{q,\nu {\rm e}^{-\tau};\mathsf{t} {\rm e}^{-\tau}}.
	\end{gather*}
	In words, the time-inhomogeneous backward $q$-Hahn process
	maps the distribution
	$\mathscr{M}^{\mathrm{qH}}_{q,\nu;\mathsf{t}}$
	onto the distribution from the same family,
	but with rescaled parameters $\mathsf{t}$ and $\nu$.
	See Fig.~$\ref{fig:two_times_qhahn_diagram}$ for an~illustration
	of the action on parameters.
\end{Theorem}

It should be noted that here
the forward $q$-Hahn TASEP has homogeneous parameters $\nu_n\equiv \nu$,
while the construction of the backward process relies on Markov swap operators
coming from inhomogeneous parameters $\nu_n$.
\begin{proof}[Proof of Theorem~\ref{thm:action_on_qhahn_distributions}]
	Take the continuous time multiparameter
	$q$-Hahn TASEP from Remark~\ref{rmk:inhom_cont_time_qhahn}
	and set the parameters to $\nu_n=\nu r^{n-1}$,
	where $\nu\in[0,1)$ and $r\in(0,1)$.
	After rescaling the continuous time by $\nu$,
	under this process each particle $x_n$ jumps by $j$ at
	rate $r^{n-1}\psi_{q,\nu r^{n-1}}(j\,|\, x_{n-1}-x_n-1)$.
	Denote the distribution of this inhomogeneous
	process at time $\mathsf{t}\in \mathbb{R}_{\ge0}$
	started from $\mathsf{step}$ by
	$\mathscr{M}^{\mathrm{qH}(r)}_{q,\nu;\mathsf{t}}$.
	Clearly,
	\begin{gather*}
		\lim_{r\to1}
		\mathscr{M}^{\mathrm{qH}(r)}_{q,\nu;\mathsf{t}}
		=
		\mathscr{M}^{\mathrm{qH}}_{q,\nu;\mathsf{t}}.
	\end{gather*}

	A suitable modification of
	Proposition~\ref{prop:qhahn_discrete_repeated_action} applies to this
	$r$-dependent distribution
	$\mathscr{M}^{\mathrm{qH}(r)}_{q,\nu;\mathsf{t}}$,
	when we take a sequence of discrete Markov backward steps:
	\begin{gather}
		\label{eq:M_r_dependent_action}
		\mathscr{M}^{\mathrm{qH}(r)}_{q,\nu;\mathsf{t}}
		\mathcal{B}^{\mathrm{qH}(r)}_{q,\nu}
		\mathcal{B}^{\mathrm{qH}(r)}_{q,r\nu}
		\mathcal{B}^{\mathrm{qH}(r)}_{q,r^2\nu}
		\cdots
		\mathcal{B}^{\mathrm{qH}(r)}_{q,r^{k-1}\nu}
		=
		\mathscr{M}^{\mathrm{qH}(r)}_{q,r^k\nu;r^k\mathsf{t}}.
	\end{gather}
	Indeed, the application of all the operators
	$\mathcal{B}^{\mathrm{qH}(r)}_{q,r^i\nu}$,
	$i=0,1,\ldots,k-1 $,
	turns
	$r^{n-1}\psi_{q,\nu r^{n-1}}$,
	the jump rate of $x_n$,
	into
	$r^{k+n-1}\psi_{q,\nu r^{k+n-1}}$.
	The latter is the same as the old jump rate
	but with the parameters $(\nu,\mathsf{t})$
	multiplied by $r^k$.
	Taking the limit as $r\to1$ in~\eqref{eq:M_r_dependent_action}
	and using Proposition~\ref{prop:limit_bwd_discrete_to_cont}
	implies the result.
\end{proof}

\subsection{Corollary. Mapping TASEP back in time}

By setting $q=\nu=0$ in Theorem~\ref{thm:action_on_qhahn_distributions},
we recover the main result of the recent paper~\cite{PetrovSaenz2019backTASEP}
on~the existence of a Markov process which maps the TASEP distributions
back in time.
Indeed, we~have for the rates~\eqref{eq:psi_things_notation}:
\begin{gather}
	\label{eq:psi_0_0}
	\psi_{0,0}(j\,|\, m)=\mathbf{1}_{j=1},\qquad
	\psi_{0,0}^{\bullet}(j'\,|\, m)=1,
\end{gather}
where $1\le j\le m$ and $0\le j'\le m-1$.
Moreover, for $q=\nu=0$
the backward continuous time process
$\mathcal{B}^{\mathrm{qH}}_{0,0}(\tau,\tau')=
\mathcal{B}^{\mathrm{qH}}_{0,0}(\tau'-\tau)$ is time-homogeneous.
Under this process, each particle~$x_n$
independently jumps to the left into one of the holes
$\{x_{n+1}+1,\ldots,x_n-1 \}$ at rate $n$ per each hole.
This dynamics is called
the \emph{backward Hammersley process}
in~\cite{PetrovSaenz2019backTASEP}.

We see from~\eqref{eq:psi_0_0} that
$\mathscr{M}^{\mathrm{qH}}_{0,0;\mathsf{t}}$
is the distribution
of the usual TASEP at time $\mathsf{t}$ started from the step initial
configuration.
Under the action of the backward Hammersley process
$\mathcal{B}^{\mathrm{qH}}_{0,0}$ for time $\tau$,
the distribution
$\mathscr{M}^{\mathrm{qH}}_{0,0;\mathsf{t}}$
maps into
$\mathscr{M}^{\mathrm{qH}}_{0,0;\mathsf{t}{\rm e}^{-\tau}}$.
This corollary of Theorem~\ref{thm:action_on_qhahn_distributions}
is precisely Theorem~1 from~\cite{PetrovSaenz2019backTASEP}.
In the latter paper the result was obtained in a completely different way
using a well-known connection (e.g., see~\cite{BorFerr2008DF})
of TASEP
and Schur processes, which are probability distributions
on two-dimensional arrays of interlacing
particles.
In contrast, here we~proved the
more general Theorem~\ref{thm:action_on_qhahn_distributions}
involving only observables of the
particle systems in~one space dimension,
and did not rely on Schur like processes in two space dimensions.

\section[Stationary dynamics on the q-TASEP distribution]
{Stationary dynamics on the $\boldsymbol q$-TASEP distribution}
\label{sec:stat_qtasep}

When $\nu=0$, the $q$-Hahn TASEP turns into the
$q$-TASEP
introduced in~\cite{BorodinCorwin2011Macdonald}.
We continue working
in the continuous time setting
as in Section~\ref{sec:qhahn_cont_time}.
In this section we introduce and study a
time-homogeneous, continuous time
Markov process which is stationary
on the distribution of~the~$q$-TASEP.

\subsection[q-TASEP and the backward process]
{$\boldsymbol q$-TASEP and the backward process}

\begin{figure}[htpb]
	\centering
	\includegraphics{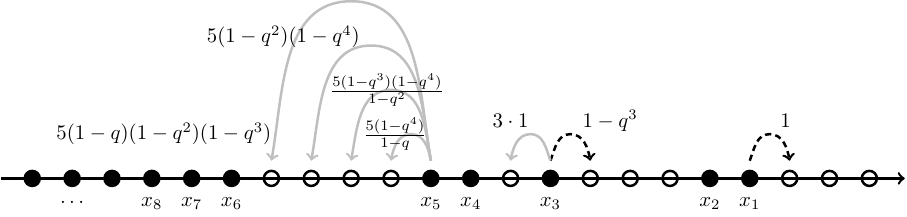}
	\caption{Examples of possible jumps with their rates in the $q$-TASEP (the dashed arrows to the right)
	and the backward $q$-TASEP (the gray arrows to the left).}
	\label{fig:qtasep_fwd_back}
\end{figure}

The \emph{$q$-TASEP} is a continuous time
Markov dynamics on $\mathrm{Conf}_{\rm fin}(\mathbb{Z})$ depending on a single parameter $q$.
Under it,
each particle $x_n(\mathsf{t})$
jumps to the right by one at rate
$1-q^{x_{n-1}(\mathsf{t})-x_n(\mathsf{t})-1}$ (by agreement, $x_0\equiv+\infty$, so the
first particle performs the Poisson random walk).
When the destination of the jump is occupied, the rate is $1-q^0=0$, so the jump is blocked
automatically.
Denote the infinitesimal Markov generator of the $q$-TASEP by
$\mathsf{T}$.
For the $q$-TASEP started from the step initial configuration
$\mathsf{step}$,
let
$\mathscr{M}^{\mathrm{qT}}_{q;\mathsf{t}}$
denote its distribution at time $\mathsf{t}\in \mathbb{R}_{\ge0}$.
See Fig.~\ref{fig:qtasep_fwd_back} (jumps to the right) for an illustration.

Setting $\nu=0$ in the backward $q$-Hahn process
from the previous
Section~\ref{sec:qhahn_cont_time}, we arrive at~the
\emph{backward $q$-TASEP}.
This
specialization in particular makes the backward process time-homogeneous.
Under the backward $q$-TASEP process (we denote its continuous time
by $\tau$), each particle $x_n(\tau)$, $n\in \mathbb{Z}_{\ge1}$,
jumps to the left into $x_n'\in \{x_{n+1}(\tau)+1,x_{n+1}(\tau)+2,\ldots,x_n(\tau)-1 \}$ at~rate
(recall notation~\eqref{eq:psi_things_notation})
\begin{gather}
	\label{eq:qtasep_back_rates}
	n\cdot \psi^{\bullet}_{q,0}(x_n'\!-\!x_{n+1}(\tau)\!-\!1\,|\, x_n(\tau)\!-\!x_{n+1}(\tau)-1)=
	\frac{n}{1\!-\!q^{x_n(\tau)-x_n'}}
	\frac{(q;q)_{x_{n}(\tau)-x_{n+1}(\tau)-1}}{(q;q)_{x_n'-x_{n+1}(\tau)-1}}.
\end{gather}
\begin{Remark}
	One can check that
	\begin{gather*}
		\sum_{x_n'=x_{n+1}+1}^{x_n-1} \psi^{\bullet}_{q,0}(x_n'-x_{n+1}-1\,|\, x_n-x_{n+1}-1)= x_n-x_{n+1}-1,
	\end{gather*}
	but we will not use this fact.
\end{Remark}
Denote the infinitesimal Markov generator of the backward $q$-TASEP by
$\mathsf{B}$.
See Fig.~\ref{fig:qtasep_fwd_back} (jumps to the left) for an illustration.

\subsection{Definition of the stationary dynamics}
\label{sub:stat_qtasep_markov_gen}

Fix the $q$-TASEP
time parameter $\mathsf{t}\in \mathbb{R}_{>0}$.
Introduce the notation
\begin{gather*}
	\mathsf{Q}^{(\mathsf{t})}:=\mathsf{T}+\frac{1}{\mathsf{t}}\,\mathsf{B}.
\end{gather*}
This is an infinitesimal Markov generator of a process
under which the
particles move to the right according to the $q$-TASEP, and
independently move to the left according to the
backward $q$-TASEP slowed down by the factor of $\mathsf{t}$.

By slightly abusing notation, we denote the continuous time Markov
process with the gene\-rator $\mathsf{Q}^{(\mathsf{t})}$
by the same letter.
We also adopt a convention
of using the letter
$\tau\in \mathbb{R}_{\ge0}$ for~the
continuous time
in the process $\mathsf{Q}^{(\mathsf{t})}$.
Thus,
$\mathsf{t}$ in $\mathsf{Q}^{(\mathsf{t})}$ is a
fixed parameter which enters the definition of the process.

\begin{Proposition}
	The continuous time Markov process with the generator $\mathsf{Q}^{(\mathsf{t})}$ is
	well-defined and can start from any configuration
	$\mathbf{x}(0)\in\mathrm{Conf}_{\rm fin}(\mathbb{Z})$.
\end{Proposition}
\begin{proof}
	The only problem in the definition of $\mathsf{Q}^{(\mathsf{t})}$
	is that it may have infinitely many jumps in finite time.
	First, observe that
	the $q$-TASEP is well-defined starting
	from any configuration $\mathbf{x}(0)\in \mathrm{Conf}_{\rm fin}(\mathbb{Z})$.
	Next, under
	$\mathsf{Q}^{(\mathsf{t})}$
	particles go to the right not faster than under the $q$-TASEP.
	Therefore, with high probability, the random configuration of particles
	under $\mathsf{Q}^{(\mathsf{t})}$
	is empty to the right and densely packed to the left
	outside a
	bounded region of $\mathbb{Z}$
	(the size of the region may depend on the time $\tau$ in $\mathsf{Q}^{(\mathsf{t})}$).
	Because of this, the total jump rate of all particles
	under~$\mathsf{Q}^{(\mathsf{t})}$
	is bounded. Therefore,
	$\mathsf{Q}^{(\mathsf{t})}$ does not
	generate infinitely many jumps in finite time
	when started from any configuration
	$\mathbf{x}(0)\in \mathrm{Conf}_{\rm fin}(\mathbb{Z})$. This completes the proof.
\end{proof}

\begin{Proposition}
	\label{prop:stationary_qtasep_process}
	The process $\mathsf{Q}^{(\mathsf{t})}$ preserves the
	$q$-TASEP distribution $\mathscr{M}_{q;\mathsf{t}}^{\mathrm{qT}}$.
\end{Proposition}
\begin{proof}
	By Theorem~\ref{thm:action_on_qhahn_distributions},
	the backward process $\mathsf{B}$, ran for small time $\delta>0$,
	maps $\mathscr{M}_{q;\mathsf{t}}^{\mathrm{qT}}$ to~$\mathscr{M}_{q;\mathsf{t}{\rm e}^{-\delta}}^{\mathrm{qT}}$.
	After evolving this distribution under the $q$-TASEP for time $\mathsf{t}-\mathsf{t}{\rm e}^{-\delta}>0$,
	we get back the original distribution
	$\mathscr{M}_{q;\mathsf{t}}^{\mathrm{qT}}$.
	Differentiating in $\delta$ and sending $\delta\searrow0$,
	the infinitesimal Markov generator of the combined dynamics
	is readily seen to be
	$\mathsf{t}\,\mathsf{T}+\mathsf{B}=\mathsf{t}\,\mathsf{Q}^{(\mathsf{t})}$.
	As the factor $\mathsf{t}$ by $\mathsf{Q}^{(\mathsf{t})}$ simply
	corresponds to the time scale (of the time variable $\tau$), we get the desired
	statement that the
	dynamics $\mathsf{Q}^{(\mathsf{t})}$
	preserves the measure $\mathscr{M}_{q;\mathsf{t}}^{\mathrm{qT}}$.
\end{proof}

\subsection[Dual process: transient q-Boson]
{Dual process: transient $\boldsymbol q$-Boson}

Our aim now is to describe the dual process
to $\mathsf{Q}^{(\mathsf{t})}$.
by means of the same duality functional
$H(\mathbf{x},\vec n)=\mathbf{1}_{n_\ell>0}\prod_{j=1}^{\ell}
q^{x_{n_j}+n_j}$
\eqref{eq:duality_functional}.
Recall that $\vec n\in \mathbb{W}^{\ell}$~\eqref{eq:W_space},
and we interpret elements of~$\mathbb{W}^{\ell}$
as $\ell$-particle configurations in $\mathbb{Z}_{\ge0}$.

\begin{figure}[htpb]
	\centering
	\includegraphics{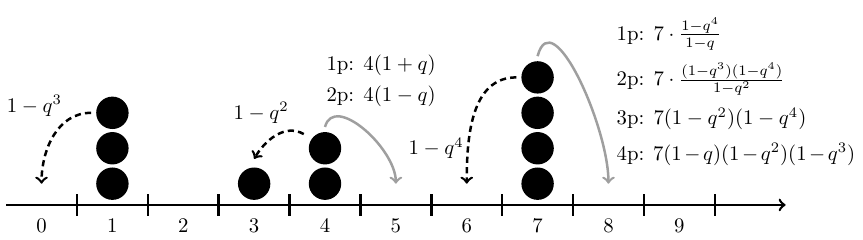}
	\caption{Examples
		of rates of possible jumps of the stochastic $q$-Boson
		$\breve{\mathsf{T}}$ (dashed jumps to the left)
		and $\breve{\mathsf{B}}$,
		the process dual to the backward $q$-TASEP (gray jumps to the right).
		Each left jump involves only one particle, while right jumps may involve any
		number of particles in the stack (in the figure,
		``1p'' means that one particle leaves the given stack, and so on).}
	\label{fig:qtasep_stationary}
\end{figure}

First consider the individual components $\mathsf{T}$ and $\mathsf{B}$ in
$\mathsf{Q}^{(\mathsf{t})}$.
The dual process to the $q$-TASEP is known as the
\emph{stochastic $q$-Boson particle system}
\cite{SasamotoWadati1998} (see also~\cite{BorodinCorwinSasamoto2012}).
Under this process, particles move in continuous time from site $k$ to $k-1$,
$k\in \mathbb{Z}_{\ge1}$,
independently
at different sites.
More precisely, if a site $k\in \mathbb{Z}_{\ge1}$ contains
$y_k$ particles, then one particle hops from site $k$ to site $k-1$
at rate $1-q^{y_k}$.
See Fig.~\ref{fig:qtasep_stationary}
(jumps to the left) for an illustration.

Denote the infinitesimal Markov generator of the stochastic $q$-Boson process by
$\breve{\mathsf{T}}$.
The duality between $\mathsf{T}$ and $\breve{\mathsf{T}}$
holds in the same sense as in Section~\ref{sub:duality}:

\begin{Proposition}[\cite{BorodinCorwinSasamoto2012}]
	We have
	\begin{gather*}
		\mathsf{T} H(\mathbf{x},\vec n)=
		\breve{\mathsf{T}} H(\mathbf{x}, \vec n)
	\end{gather*}
	for any $\mathbf{x}\in \mathrm{Conf}_{\rm fin}(\mathbb{Z})$ and $\vec n\in \mathbb{W}^{\ell}$,
	where $\mathsf{T}$ acts on $\mathbf{x}$, and $\breve{\mathsf{T}}$ acts on $\vec n$.
\end{Proposition}
Let us now define another continuous time Markov process on $\mathbb{W}^{\ell}$
which will be dual to the backward $q$-TASEP.
Under this new process, particles move in continuous time from site $k$ to~$k+1$,
$k\in \mathbb{Z}_{\ge1}$,
independently
at different sites.
More precisely, if a site $k\in \mathbb{Z}_{\ge0}$ contains
$y_k$ particles, then the process sends $y_k-j$ particles to site $k+1$,
where $j\in \{0,1,\ldots,y_{k-1} \}$,
at rate $k\cdot \psi^{\bullet}_{q,0}(j\,|\, y_k)$.
(In particular, particles cannot leave site $0$.)
Denote the infinitesimal Markov generator of this dynamics
by $\breve{\mathsf{B}}$.
See Fig.~\ref{fig:qtasep_stationary} (right jumps)
for an illustration.

\begin{Proposition}
	The Markov generator $\breve{\mathsf{B}}$
	is dual to the backward $q$-TASEP generator $\mathsf{B}$:
	\begin{gather*}
		\mathsf{B} H(\mathbf{x},\vec n)=
		\breve{\mathsf{B}} H(\mathbf{x}, \vec n)
	\end{gather*}
	for any
	$\mathbf{x}\in \mathrm{Conf}_{\rm fin}(\mathbb{Z})$ and $\vec n\in \mathbb{W}^{\ell}$,
	where $\mathsf{B}$ acts on $\mathbf{x}$, and $\breve{\mathsf{B}}$ acts on $\vec n$.
	Consequently, the infinitesimal Markov generator
	\begin{gather*}
		\breve{\mathsf{Q}}^{(\mathsf{t})}:=
		\breve{\mathsf{T}}+\frac{1}{\mathsf{t}}\breve{\mathsf{B}}
	\end{gather*}
	is dual to $\mathsf{Q}^{(\mathsf{t})}$, the generator of the stationary dynamics
	on the $q$-TASEP distribution.
\end{Proposition}
\begin{proof}
	The duality between $\mathsf{B}$ and $\breve{\mathsf{B}}$
	follows from Proposition~\ref{prop:qhahn_transition_duality}
	via the continuous time limit described
	in Section~\ref{sec:qhahn_cont_time}.
	Indeed, the backward $q$-Hahn process
	is
	a continuous time limit of the
	combination
	of the steps $p_k^{\mathrm{qH}}$
	applied at each site of the lattice $\mathbb{Z}$.
	The dual process~$\breve{\mathsf{B}}$ is the limit of the same type of~the
	combination of the dual steps
	$\breve p_k^{\mathrm{qH}}$, together with the dege\-ne\-ra\-tion~$\nu=0$.
	This implies the duality of
	$\mathsf{B}$ and $\breve{\mathsf{B}}$.
	The claim about the duality of
	$\mathsf{Q}^{(\mathsf{t})}$
	and
	$\breve{\mathsf{Q}}^{(\mathsf{t})}$
	follows by linearity.
\end{proof}

Because the rates of the right jumps
under $\breve{\mathsf{Q}}^{(\mathsf{t})}$
grow as the particles of $\vec n$ get farther to the right, the
process $\breve{\mathsf{Q}}^{(\mathsf{t})}$ is transient
(except the absorption at $n_\ell=0$),
see the proof of Lemma~\ref{lemma:limit_of_moments}
below for details.
For this reason we call $\breve{\mathsf{Q}}^{(\mathsf{t})}$
the \emph{transient stochastic $q$-Boson} particle system on $\mathbb{W}^{\ell}$ (transient
$q$-Boson, for short).
Let us record the duality between it and $\mathsf{Q}^{(\mathsf{t})}$
in terms of~expectations:
\begin{Corollary}
	\label{cor:duality_Q_expectations}
	Fix $\mathsf{t}\in \mathbb{R}_{>0}$ and
	$\ell\in \mathbb{Z}_{\ge1}$.
	Take
	any $\mathbf{x}^{0}\in \mathrm{Conf}_{\rm fin}(\mathbb{Z})$,
	and any
	$\vec n^0\in \mathbb{W}^{\ell}$.
	Let $\{\mathbf{x}(\tau)\}_{\tau\in \mathbb{R}_{\ge0}}$ be the
	process on $\mathrm{Conf}_{\rm fin}(\mathbb{Z})$
	with generator $\mathsf{Q}^{(\mathsf{t})}$
	started from $\mathbf{x}^{0}$, and
	$\{\vec n(\tau)\}_{\tau\in \mathbb{R}_{\ge0}}$
	be the process on $\mathbb{W}^{\ell}$
	with generator $\breve{\mathsf{Q}}^{(\mathsf{t})}$
	started from $\vec n^0$.
	Then for any $\tau\in \mathbb{R}_{\ge0}$ we have
	\begin{gather}
		\mathbb{E}^{\mathrm{stat}(\mathsf{t})}_{\mathbf{x}(0)=\mathbf{x}^0}
		H(\mathbf{x}(\tau),\vec n^0)=
		\mathbb{E}^{\mathrm{trqBoson}(\mathsf{t})}_{\vec n(0)=\vec n^0}
		H(\mathbf{x}^0,\vec n(\tau)).
	\end{gather}
\end{Corollary}
The example of this duality statement (and further discussion
towards the results of the next
Section~\ref{sub:stat_qtasep_convergence})
in the simplest case $\ell=1$ may be found in
Section~\ref{sub:intro_stat_tasep_briefly}
in Introduction.

\subsection{Convergence to the stationary distribution}
\label{sub:stat_qtasep_convergence}

In this section we use duality to prove the following
result:
\begin{Theorem}
	\label{thm:qtasep_converges_to_stationary}
	Fix $\mathsf{t}\in \mathbb{R}_{>0}$.
	For any initial configuration $\mathbf{x}(0)\in \mathrm{Conf}_{\rm fin}(\mathbb{Z})$
	the Markov process~$\mathbf{x}(\tau)$ with the generator
	$\mathsf{Q}^{(\mathsf{t})}$
	converges, as $\tau\to+\infty$,
	to the stationary distribution $\mathscr{M}^{\mathrm{qT}}_{q;\mathsf{t}}$
	$($in~the sense of joint distributions of arbitrary
	finite subcollections of particles$)$.
\end{Theorem}
The proof of Theorem~\ref{thm:qtasep_converges_to_stationary}
occupies the rest of this section.

Fix any $\ell\ge1$ and $\vec m\in \mathbb{W}^{\ell}$.
Let $\vec n(\tau)$ be the
transient $q$-Boson
$\breve{\mathsf{Q}}^{(\mathsf{t})}$ started from $\vec m$.
Denote the \emph{survival probability} of the transient $q$-Boson till time $\tau$ by
\begin{gather*}
	S_\tau(\vec m):=\mathbb{P}(n_\ell(\tau)>0\,|\, \vec n(0)=\vec m).
\end{gather*}
Note that if $m_\ell=0$, we automatically have $S_\tau(\vec m)=0$ for all $\tau$.

\begin{Lemma}
	For any $\vec m\in \mathbb{W}^{\ell}$,
	the asymptotic survival probability
	\begin{gather*}
		S(\vec m):=\lim\limits_{\tau\to+\infty}S_\tau(\vec m)
	\end{gather*}
	exists.
\end{Lemma}
\begin{proof}
	We have
	\begin{gather*}
		S_\tau(\vec m)=\mathbb{P}(n_\ell(\tau)>0\,|\, \vec n(0)=\vec m)=
		\mathbb{P}(n_\ell(s)>0\textnormal{ for all $s\in [0,\tau]$}\,|\, \vec n(0)=\vec m)
	\end{gather*}
	because once $n_\ell$ reaches $0$, it can never become positive again.
	Therefore, the quantities~$S_{\tau}(\vec m)$ decrease in $\tau$ due to monotonicity
	in $\tau$
	of the
	events
	$\{n_\ell(s)>0\textnormal{ for all $s\in [0,\tau]$}\}$.
	Because $S_\tau(\vec m)$'s are nonnegative, the limit exists.
\end{proof}

Observe that $H(\mathsf{step},\vec m)=\mathbf{1}_{m_\ell>0}$.
Therefore, for the stationary process
$\mathbf{x}(\tau)$ started from $\mathsf{step}$ we have
$\mathbb{E}^{\mathrm{stat}(\mathsf{t})}_{\mathsf{step}}
\prod_{j=1}^{\ell}q^{x_{m_j}(\tau)+m_j}=S_\tau(\vec m)$
for all $\vec m\in \mathbb{W}^{\ell}$
(in particular, if $m_\ell=0$ then both sides are zero).
For other initial conditions for $\mathsf{Q}^{(\mathsf{t})}$
this identity does not hold for finite time $\tau$, but it still
holds asymptotically:

\begin{Lemma}
	\label{lemma:limit_of_moments}
	For any $\vec m\in \mathbb{W}^\ell$
	and any initial data $\mathbf{x}^0\in \mathrm{Conf}_{\rm fin}(\mathbb{Z})$
	for the stationary process~$\mathsf{Q}^{(\mathsf{t})}$,
	we have
	\begin{gather}
		\label{eq:qmoments_converge_to_survival_probability}
		\lim_{\tau\to+\infty}
		\mathbb{E}_{\mathbf{x}(0)=\mathbf{x}^0}^{\mathrm{stat}(\mathsf{t})}
		\prod_{j=1}^{\ell}
		q^{x_{m_j}(\tau)+m_j}=S(\vec m),
	\end{gather}
	where $S(\vec m)$ is the asymptotic survival probability
	of the transient $q$-Boson.
\end{Lemma}
\begin{proof}
	If $m_\ell=0$, then both sides of
	\eqref{eq:qmoments_converge_to_survival_probability}
	are zero, so it suffices to assume
	that $m_\ell>0$.
	By~dua\-lity (Corollary~\ref{cor:duality_Q_expectations}),
	the left-hand side of
	\eqref{eq:qmoments_converge_to_survival_probability}
	is equal to
	\begin{gather}
		\label{eq:limit_of_moments_proof}
		\lim_{\tau\to+\infty}
		\mathbb{E}^{\mathrm{trqBoson}(\mathsf{t})}
		_{\vec n(0)=\vec m}
		\bigg(\mathbf{1}_{n_\ell(\tau)>0}\prod_{j=1}^{\ell}
		q^{x_{n_j(\tau)}^0+n_j(\tau)}\bigg).
	\end{gather}
	Here $x_{n}^{0}$ are the particle
	coordinates under the initial data $\mathbf{x}^{0}$.

	First, we show that the Markov process
	$\vec n(\tau)$ conditioned
	to stay in the region $\left\{ n_\ell\ge1 \right\}$
	is transient.
	Observe that
	we can couple the $\ell$-particle
	process
	$\vec n(\tau)$ restricted to this region
	with a single-particle process
	$Y(\tau)$
	on $\mathbb{Z}_{\ge1}$
	with jump rates
	\begin{gather*}
		\mathrm{rate}_Y(k+1\to k)=1-q^{\ell},\quad
		\mathrm{rate}_Y(k\to k+1)=\frac{k}{\mathsf{t}}
		\Big(\min_{j,r\colon 0\le j<r\le \ell}
		\psi^{\bullet}_{q,0}(j\,|\, r)\Big)>0,\quad
		k\ge1.
	\end{gather*}
	The coupling is such that $Y(0)=n_\ell(0)$,
	all left jumps of $n_\ell$ force $Y$ to jump to the left,
	and all right jumps of $Y$ force $n_\ell$ to jump to the right.
	This implies that under this coupling
	we have $Y(\tau)\le n_\ell(\tau)$ for all $\tau>0$.

	The process $Y(\tau)$ on $\mathbb{Z}_{\ge1}$
	is a standard example of a transient Markov process:
	it eventually (with probability $1$)
	reaches the
	part $\{A,A+1,\ldots \}\subset\mathbb{Z}_{\ge1}$ (where $A$
	depends only
	on $\ell$ and $\mathsf{t}$)
	where the average drift to the right is bounded away from zero.
	With positive probability $Y(\tau)$ then never comes back
	from $\left\{ A,A+1,\ldots \right\}$ to the neighborhood of zero,
	and thus is transient.

	Using the transience, we can lower bound $S(\vec m)$
	by the (positive) probability of the event that:
	(1) if there were any particles at site $1$ at time $\tau=0$,
	then all these particles leave $1$ by~right jumps;
	(2) after that, $n_\ell(\tau)$
	never comes back to $1$ (and hence all other particles of $\vec n(\tau)$ also never come
	back to $1$).
	This implies that $S(\vec m)>0$.

Now denote by $\mathscr{S}(\vec m)$ the event that $n_\ell(\tau)>0$ for all $\tau$
	conditioned on $\vec n(0)=\vec m$. Thus,
	$\mathbb{P}(\mathscr{S}(\vec m))=S(\vec m)>0$.
	We have for the expectation in~\eqref{eq:limit_of_moments_proof}:
	\begin{gather*}
		\mathbb{E}^{\mathrm{trqBoson}(\mathsf{t})}		_{\vec n(0)=\vec m}
		\bigg(\mathbf{1}_{n_\ell(\tau)>0}\prod_{j=1}^{\ell}
		q^{x_{n_j(\tau)}^0+n_j(\tau)}\bigg)		=
		S(\vec m)\,		\mathbb{E}^{\mathrm{trqBoson}(\mathsf{t})}
		_{\vec n(0)=\vec m}		\bigg(\prod_{j=1}^{\ell}
		q^{x_{n_j(\tau)}^0+n_j(\tau)}\,|\, \mathscr{S}(\vec m)\bigg)
		\\ \hphantom{\mathbb{E}^{\mathrm{trqBoson}(\mathsf{t})}		_{\vec n(0)=\vec m}}
{}+( 1-S(\vec m) )\mathbb{E}^{\mathrm{trqBoson}(\mathsf{t})}_{\vec n(0)=\vec m}
		\bigg(\mathbf{1}_{n_\ell(\tau)>0}\prod_{j=1}^{\ell}
		q^{x_{n_j(\tau)}^0+n_j(\tau)}\,|\, \mathscr{S}(\vec m)^{c}\bigg).
	\end{gather*}
	The second summand goes to zero as $\tau\to+\infty$,
	since inside the event $\mathscr{S}(\vec m)^c$, we have
	$n_\ell(s)=0$ for all $s\in(s_0,+\infty)$ (where $s_0$ is random but finite).
	In the first summand, observe that
	conditioned on the asymptotic survival $\mathscr{S}(\vec m)$,
	we have almost surely due to transience that
	$n_j(\tau)\to +\infty$, $\tau\to+\infty$,
	for all
	$j=1,\ldots,\ell $.
	Because the initial configuration $\mathbf{x}^0\in \mathrm{Conf}_{\rm fin}(\mathbb{Z})$
	is densely packed to the left, we thus
	almost surely have
	$q^{x_{n_j(\tau)}^0+n_j(\tau)}\to 1$ for all $j$ as $\tau\to+\infty$.
	Therefore, the expectation of the product in the first summand
	tends to $1$, and we see that~\eqref{eq:limit_of_moments_proof} is equal
	to~$S(\vec m)$, as desired.
\end{proof}

The asymptotic survival probabilities $S(\vec m)$
satisfy certain \emph{normalization at infinity}:
\begin{Lemma}
	\label{lemma:S_m_normalization_at_infinity}
	Fix $\ell$. For any $\varepsilon>0$ there exists
	$R=R(\ell,\varepsilon)\in \mathbb{Z}_{\ge1}$ such that
	for all $\vec m\in \mathbb{W}^{\ell}$
	with $m_1>R$ we have
	\begin{gather*}
		|S(m_1,m_2,\ldots,m_\ell )-S(m_2,\ldots,m_\ell )|<\varepsilon,
	\end{gather*}
	where
	$S(m_2,\ldots,m_\ell )$ is the survival probability of the
	transient $q$-Boson on $\mathbb{W}^{\ell-1}$.
	If $\ell=1$, then $S(m_2,\ldots,m_\ell )=1$ by agreement.
\end{Lemma}
\begin{proof}
	We can assume that $m_\ell>0$, otherwise both expressions
	$S(\cdot)$ in the claim are zero.
	The desired statement follows from the transience
	of the process $\vec n(\tau)$ (with generator $\breve{\mathsf{Q}}^{(\mathsf{t})}$)
	as in the proof of the previous Lemma~\ref{lemma:limit_of_moments}.
	Namely, we can lower bound the jump rate of $n_1(\tau)$ to the right
	from a site $k\in \mathbb{Z}_{\ge1}$ by
	$\mathrm{const}\cdot k$.
	Let $\vec n(\tau)$ start from $\vec m$.
	If $m_1>R$ is large,
	the probability that the first particle
	$n_1(\tau)$ ever returns to the $R/2$-neighborhood of zero is
	close to zero. Thus,
	the
	probability $S(\vec m)$
	that the process $\vec n(\tau)$ started from $\vec m$
	survives and runs off to infinity
	is close to the asymptotic survival probability
	of the process on $\mathbb{W}^{\ell-1}$ with one less particle and started from
	$(m_2,\ldots,m_\ell )$.
	In the special case $\ell=1$,
	the claim reads
	$S(m_1)\to 1$
	as $m_1\to+\infty$, which clearly holds.
	This implies the claim.
\end{proof}
\begin{Remark}
	If under the conditions of Lemma~\ref{lemma:S_m_normalization_at_infinity}
	the second coordinate $m_2$ is also very large, then
	one can similarly show that
	both $S(m_1,m_2,m_3,\ldots,m_\ell )$ and
	$S(m_2,m_3,\ldots,m_\ell )$
	are close to
	$S(m_3,\ldots,m_\ell )$,
	and thus close to each other.
	Hence an analogue of
	Lemma~\ref{lemma:S_m_normalization_at_infinity}
	for a~number of first coordinates $m_1,\ldots,m_j $
	being large also holds.
\end{Remark}

To finish the proof of
Theorem~\ref{thm:qtasep_converges_to_stationary}
it remains to show that
\begin{gather}
	\label{eq:qtasep_qmoments_notation_in_stationary_stuff}
	S(\vec m)=
	\mathbb{E}^{\mathrm{qT}}_{\mathsf{step}}\prod_{j=1}^{\ell}q^{x_{m_j}(\mathsf{t})+m_j}
\end{gather}
for all $\ell$ and $\vec m\in \mathbb{W}^{\ell}$.
Here the quantity in the left-hand side is the
long time limit of the $q$-moment of
$\mathsf{Q}^{(\mathsf{t})}$, and
the right-hand side is the $q$-moment of the
$q$-TASEP. This suffices since
in our situation the $q$-moments uniquely
characterize the distribution.
We will establish~\eqref{eq:qtasep_qmoments_notation_in_stationary_stuff}
by~showing that both sides satisfy the same
equations (harmonicity with respect to the transient $q$-Boson)
plus normalization at infinity which
uniquely determine the function.

\begin{Lemma}
	As function of $\vec m$, the survival probabilities $S(\vec m)$
	are harmonic for the transient $q$-Boson, that is,
	\begin{gather}
		\label{eq:S_m_action_annihilated}
		(\breve{\mathsf{Q}}^{(\mathsf{t})}S)(\vec m)=
		\sum_{\vec m'}
		\breve{\mathsf{Q}}^{(\mathsf{t})}(\vec m,\vec m')
		S(\vec m')=0\qquad
		\textnormal{for all $\vec m\in \mathbb{W}^{\ell}$}.
	\end{gather}
	Here the first identity is simply the expression for the action of the
	generator on a function, and the claim is that this action gives identical $0$.
\end{Lemma}
\begin{proof}
	The argument is rather standard.
	Consider the evolution of the process $\vec n(\cdot)$
	started from $\vec m$ during short time ${\rm d}\tau$.
	Conditioned that the process stepped into $\vec m'$ (which happens with probability
	$\breve{\mathsf{Q}}^{(\mathsf{t})}(\vec m,\vec m'){\rm d}\tau$),
	the survival probability is then $S(\vec m')$.
	With the complementary probability
	$1-\sum_{\vec m'\colon \vec m'\ne \vec m}\breve{\mathsf{Q}}^{(\mathsf{t})}(\vec m,\vec m'){\rm d}\tau=
	1+\breve{\mathsf{Q}}^{(\mathsf{t})}(\vec m,\vec m){\rm d}\tau$, the
	process did not leave $\vec m$, and the survival probability did not change.
	Therefore,
	\begin{gather*}
		S(\vec m)=S(\vec m)+
		{\rm d}\tau
		\sum_{\vec m'}
		\breve{\mathsf{Q}}^{(\mathsf{t})}(\vec m,\vec m')
		S(\vec m').
	\end{gather*}
	Taking the coefficient by ${\rm d}\tau$
	leads to the desired identity.
\end{proof}

\begin{Lemma}
	\label{lemma:Sm_equations_have_unique_solution}
	Harmonicity condition~\eqref{eq:S_m_action_annihilated}
	together with
	Lemma~$\ref{lemma:S_m_normalization_at_infinity}$
	$($normalization at infi\-nity$)$
	and the condition that $S(\vec m)=0$ whenever $m_\ell=0$
	uniquely determine the function
	$S(\vec m)$, $\vec m\in \mathbb{W}^{\ell}$.
\end{Lemma}
\begin{proof}
	Assume that $G_\ell(\vec m)$,
	where $\ell=1,2,\ldots $ and
	$\vec m\in \mathbb{W}^{\ell}$, is a
	family of harmonic functions
	satisfying normalization at infinity as
	in Lemma~\ref{lemma:S_m_normalization_at_infinity}, that is,
	for any $\varepsilon>0$ there exists $R$ such that
	for all $\vec m\in \mathbb{W}^{\ell}$ with $m_1>R$ we have
	\begin{gather}
		\label{eq:last_lemma_proof05}
		|G_\ell(m_1,m_2,\ldots,m_\ell )-G_{\ell-1}(m_2,\ldots,m_\ell )|<\varepsilon.
	\end{gather}
	Moreover, we assume that $G_\ell(\vec m)=0$ whenever $m_\ell=0$.
	We will argue by induction on $\ell$
	and show that $G_{\ell}(\vec m)$ is equal to $S(\vec m)$, the
	asymptotic survival probability of the transient
	$q$-Boson on $\mathbb{W}^{\ell}$.

	For $\ell=1$, it is straightforward to see
	that the
	space of harmonic functions
	vanishing at $0$
	is one-dimensional.
	In this case the normalization at infinity is the single condition
	$G_1(m_1)\to 1$ as $m_1\to+\infty$,
	which determines the harmonic function uniquely.

	Next, observe that
	because the function $G_\ell$ is harmonic, it satisfies
	the following averaging property:
	\begin{gather*}
		G_\ell(\vec m)=\mathbb{E}^{\mathrm{trqBoson}(\mathsf{t})}_{\vec n(0)=\vec m}
		G_\ell\left( \vec n(\tau) \right),
	\end{gather*}
	where $\tau$ is arbitrary.

	Assume that $m_\ell>0$ (this does not restrict the generality).
	For $R\in \mathbb{Z}_{\ge1}$, take the stopping time
	\begin{gather*}
		T_R:=\inf\{\tau\ge0 \colon n_\ell(\tau)=0\ \textnormal{or}
		\ n_1(\tau)=R\},
	\end{gather*}
	where the process $\vec n(\tau)$ starts from $\vec m$.
	This stopping time is almost surely bounded and has finite expectation.
	Then
	\begin{gather}
		\label{eq:last_lemma_proof}
		G_\ell(\vec m)=
		\mathbb{E}^{\mathrm{trqBoson}(\mathsf{t})}_{\vec n(0)=\vec m}
		G_\ell(\vec n(T_R) )=
		\mathbb{E}^{\mathrm{trqBoson}(\mathsf{t})}_{\vec n(0)=\vec m}
		\{G_\ell(R,n_2(T_R),\ldots, n_\ell(T_R))\mathbf{1}_{n_1(T_R)=R}\}.
	\end{gather}
	The first equality above follows from
	the optional stopping theorem, and the second
	one is the splitting into two cases, $n_1(T_R)=R$ or $n_\ell(T_R)=0$.
	In the latter case the function $G_\ell$ vanishes by our assumptions.

	Take any $\varepsilon>0$ and choose $R$ such that~\eqref{eq:last_lemma_proof05}
	holds.
	Thus, by the normalization at infinity, we have
	\begin{gather*}
		\left|{\rm RHS}\eqref{eq:last_lemma_proof}-
		\mathbb{E}^{\mathrm{trqBoson}(\mathsf{t})}_{\vec n(0)=\vec m}
		\{S(n_2(T_R),\ldots, n_\ell(T_R))\mathbf{1}_{n_1(T_R)=R}\}\right|<\varepsilon,
	\end{gather*}
	where we have replaced $G_{\ell-1}$ by $S$ using the induction hypothesis.
	Since $S$ satisfies the same conditions
	(harmonicity, normalization at infinity, and vanishing when $m_\ell=0$)
	as the family~$G_\ell$, identity~\eqref{eq:last_lemma_proof}
	is also valid for $S$. Thus,
	\begin{gather*}
		\left|\mathbb{E}^{\mathrm{trqBoson}(\mathsf{t})}_{\vec n(0)=\vec m}
		\{S(n_2(T_R),\ldots, n_\ell(T_R))\mathbf{1}_{n_1(T_R)=R}\}
		-S(m_1,\ldots,m_\ell)\right|<\varepsilon,
	\end{gather*}
	which means that
	$G_\ell(m_1,\ldots,m_\ell )$ is $\varepsilon$-close to $S(m_1,\ldots,m_\ell )$.
	This completes the proof.
\end{proof}

\begin{Lemma}
	The $q$-moments of the $q$-TASEP in the right-hand side of
	\eqref{eq:qtasep_qmoments_notation_in_stationary_stuff}
	satisfy
	all the conditions listed in the previous Lemma~$\ref{lemma:Sm_equations_have_unique_solution}$.
\end{Lemma}
\begin{proof}
	The harmonicity of the $q$-moments
	follows from the fact that the $q$-TASEP distribution~$\mathscr{M}^{\mathrm{qT}}_{q;\mathsf{t}}$
	is stationary under $\mathsf{Q}^{(\mathsf{t})}$
	(Proposition~\ref{prop:stationary_qtasep_process}), together with duality
	between $\mathsf{Q}^{(\mathsf{t})}$ and the transient
	$q$-Boson
	(Corollary~\ref{cor:duality_Q_expectations}).

	The normalization at infinity follows from the fact that the
	$q$-TASEP started from $\mathsf{step}$
	lives on $\mathrm{Conf}_{\rm fin}(\mathbb{Z})$,
	so for any $\mathsf{t}\in \mathbb{R}_{\ge0}$ we almost surely have
	$x_{m}(\mathsf{t})+m\to 0$ as $m\to+\infty$.

The fact that the $q$-moments vanish when $n_\ell=0$
	follows from the agreement that \linebreak \mbox{$x_0=+\infty$}.
\end{proof}

Combining all the lemmas in this section,
we get the desired
Theorem~\ref{thm:qtasep_converges_to_stationary}.

\section{Beta polymer}
\label{sec:beta_polymer}

The $q$-Hahn TASEP has a remarkable degeneration~-- the beta polymer model introduced in~\cite{CorwinBarraquand2015Beta}.
This model
is also related to a random walk in dynamic beta random environment,
but here we will formulate everything only in terms of the polymer model.
In this section we
present Markov swap operators for the
multiparameter beta polymer.
The swap operators can be realized as certain additional
layers in the strict-weak lattice on which the beta polymer
is defined.

\subsection{Multiparameter beta polymer and its joint moments}

\begin{figure}[htpb]
	\centering
	\includegraphics{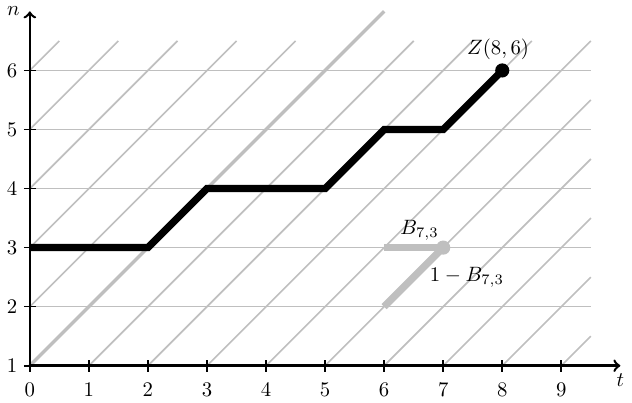}
	\caption{The beta polymer partition function
		$Z(t,n)$ is the sum of weights of all
		strict-weak paths from the line $\{0\}\times \mathbb{Z}_{\ge1}$ to $(t,n)$, where the
		weights of the horizontal and the diagonal edges are
		$B_{\cdots}$ and $1-B_{\cdots}$, respectively.
		The weight of a path is defined as the product of its edge weights.
		The initial condition along the left boundary is $Z(0,n)=1$ for all $n$,
		so $Z(t,n)=1$ above the diagonal (i.e.,~for~$n>t$).}
	\label{fig:beta_polymer}
\end{figure}

Take parameters $\upgamma>0$ and $\upnu_n>0$, $t,n\in \mathbb{Z}_{\ge1}$,
such that
\begin{gather*}
	\textnormal{$\upnu_n-\upgamma>0$\quad
	for all $n$},\qquad
	\textnormal{$\upnu_i-\upnu_j\notin\mathbb{Z}$\quad for all $i$, $j$}.
\end{gather*}
Let $B_{t,n}\sim \mathrm{Beta}(\upnu_n-\upgamma,\upgamma)$ be independent
beta distributed random variables.
Here by the beta distribution $\mathrm{Beta}(\alpha,\beta)$ we mean
the one with the density
\begin{gather*}
	\frac{\Gamma(\alpha+\beta)}{\Gamma(\alpha)\Gamma(\beta)}\,x^{\alpha-1}(1-x)^{\beta-1},
	\qquad
	x\in [ 0,1 ].
\end{gather*}
The (inhomogeneous) beta polymer $\{Z(t,n)\}_{t,n\in \mathbb{Z}_{\ge1}}$
is a collection of random variables
satisfying the random recursion
\begin{gather*}
	\begin{split}
		Z(t,n)&=B_{t,n}Z(t-1,n)+(1-B_{t,n})Z(t-1,n-1),\\
		Z(t,1)&=B_{t,1}Z(t-1,1),
	\end{split}
\end{gather*}
with the initial condition
$Z(0,n)=1$ for all $n\in \mathbb{Z}_{\ge1}$.
See Fig.~\ref{fig:beta_polymer} for a graphical interpretation of the beta
polymer as a point-to-line partition function on the
strict-weak lattice.

The homogeneous version of the beta polymer
was studied in~\cite{CorwinBarraquand2015Beta}
as a scaling limit $q\to1$ of~the $q$-Hahn TASEP.
The multiparameter generalization
is obtained through the same limit from our model described in Section~\ref{sec:qhahn_TASEP}.

The random variables $Z(t,n)$ are between
$0$ and $1$.\footnote{Moreover,
one can check that for fixed $t$ they are ordered as
$1\ge Z(t,1)\ge Z(t,2)\ge Z(t,3)\ge \cdots $.}
Because of this, the joint distribution of~the beta polymer random variables $\{Z(t,n)\}_{n\in \mathbb{Z}_{\ge1}}$
(for every fixed $t$)
is determined by the joint moments.
These moments have the following form:
\begin{Proposition}
	For each $t\in \mathbb{Z}_{\ge1}$
	and any $n_1\ge \dots\ge n_k\ge1 $ we have
	\begin{gather}
	\mathbb{E}^{\mathrm{beta}(\boldsymbol\upnu)}(Z(t,n_1)\cdots Z(t,n_k))\nonumber
	\\ \qquad
 {}=\frac{1}{(2\pi i)^k}\oint\dots\oint
	\prod_{1\le A<B\le k}\frac{z_A-z_B}{z_A-z_B-1}\prod_{j=1}^{k}
	\bigg(\prod_{i=1}^{n_j}\frac{z_j}{z_j-\upnu_i}\bigg)
	\bigg(\frac{z_j-\upgamma}{z_j} \bigg)^t\frac{{\rm d}z_j}{z_j}.
	\label{eq:beta_polymer_moments}
	\end{gather}
	The integration contours
	are around $\{\upnu_i\}$,
	do not encircle $0$, and
	the contour for $z_j$ contains the contour for $z_{j+1}+1$,
	$j=1,\ldots,k-1 $.
\end{Proposition}
\begin{proof}
	For $\upnu_i\equiv \upnu$, this formula is a simple change of variables
	from
	\cite[Proposition 2.11]{CorwinBarraquand2015Beta}, where $\upmu=\upnu-\upgamma$.
	Its generalization with different $\upnu_i$'s
	is obtained by duality and coordinate Bethe ansatz in the same manner
	as in~\cite{CorwinBarraquand2015Beta}
	by checking that the
	right-hand side of~\eqref{eq:beta_polymer_moments}
	satisfies the free evolution equations, boundary conditions, and
	the initial condition at $t=0$.
\end{proof}

For future use, let us recall here the free evolution
equations and the boundary conditions satisfied by the
right-hand side of~\eqref{eq:beta_polymer_moments}.
Denote this right-hand side by $f(t;n_1,\ldots,n_k )$.
Then it satisfies the two-body boundary conditions
\begin{gather}
	f(t;n_1,\ldots,n_{i}-1,n_{i+1}-1,\ldots,n_k )+
	(\upnu_{n_i}-1)f(t;n_1,\ldots,n_{i},n_{i+1}-1,\ldots,n_k )\nonumber
	\\ \qquad
	{}+	f(t;n_1,\ldots,n_{i},n_{i+1},\ldots,n_k )-
	(\upnu_{n_i}+1)f(t;n_1,\ldots,n_{i}-1,n_{i+1},\ldots,n_k )=0
		\label{eq:beta_boundary_conditions}
\end{gather}
for all $\vec n\in \mathbb{Z}^{k}$
such that for some $i\in \left\{ 1,\ldots,k \right\}$,
$n_i=n_{i+1}$.
This is checked similarly to the homogeneous case \cite[Section 4.1]{CorwinBarraquand2015Beta},
because the inhomogeneity parameters $\upnu_{n_i}$ are the same for~each cluster of equal $n_j$'s
(a similar effect is observed in the proof of Lemma~\ref{lemma:boundary_conditions} where
only one of the parameters, $\nu_n$, plays an essential role).

The free evolution equations satisfied by $f(t;\vec n)$ are
\begin{gather*}
	f(t+1;\vec n)=
	\prod_{i=1}^{k}\big[\nabla^{\textnormal{beta}}_{\upgamma/\upnu_{n_i}}\big]_{i}\,
	f(t;\vec n),
\end{gather*}
where
$\nabla^{\textnormal{beta}}_p g(n):=p \, g(n-1)+(1-p)\,g(n)$,
and the operator $[\nabla^{\textnormal{beta}}_{\upgamma/\upnu_{n_i}}]_{i}$
is applied in the $i$-th variable $n_i$.

\begin{Remark}
	Similarly to Remark~\ref{rmk:inhom_time},
	one can generalize the
	beta polymer model so that the
	parameter $\upgamma$ depends on $t$.
	The moment formula~\eqref{eq:beta_polymer_moments}
	and our main result (Theorem~\ref{thm:beta_polymer_back} below) would continue to hold
	with straightforward modifications. For simplicity, everywhere below
	we take $\upgamma$ independent of $t$.
\end{Remark}

\begin{Proposition}
	The multiparameter beta polymer,
	viewed as a stochastic particle system
	$Z(t,\cdot)$ on $[0,1]$ with time $t\in \mathbb{Z}_{\ge1}$,
	is parameter-symmetric in the sense of Definition~$\ref{def:symm_IPS}$.
	Moreover, the distribution of each $Z(t,n)$
	depends on the parameters $\upnu_1,\ldots,\upnu_n $ in a symmetric way.
\end{Proposition}

\subsection{Swap operator for the beta polymer}

Let $\{Z^{\boldsymbol\upnu}(t,n)\}_{t,n\in \mathbb{Z}_{\ge1}}$
be the beta polymer
with parameters $\boldsymbol\upnu=\{\upnu_1,\upnu_2,\ldots \}$.
Recall that by~$s_n$ we denote the elementary transposition
$(n,n+1)$. The next statement presents
the swap operator interchanging
$\upnu_n \leftrightarrow\upnu_{n+1}$.

\begin{Theorem}
	\label{thm:beta_polymer_back}
	Fix $t,n\in \mathbb{Z}_{\ge1}$
	and assume that $\upnu_n<\upnu_{n+1}$.
	Let $\tilde B\sim \mathrm{Beta}(\upnu_{n+1}-\upnu_n,\upnu_n)$ be a new
	beta random variable independent
	of the environment $\{B_{t,n}\}$ (and hence of the beta polymer).
	Then we have equality in distribution
	\begin{multline*}
		( Z^{\boldsymbol\upnu}(t,1),\ldots,
			Z^{\boldsymbol\upnu}(t,n-1),
			\tilde B Z^{\boldsymbol\upnu}(t,n+1)+(1-\tilde B)Z^{\boldsymbol\upnu}(t,n),
		Z^{\boldsymbol\upnu}(t,n+1))\\\stackrel{d}{=}
		(
			Z^{s_n \boldsymbol\upnu}(t,1),
			\ldots,
			Z^{s_n \boldsymbol\upnu}(t,n-1),
			Z^{s_n \boldsymbol\upnu}(t,n),
			Z^{s_n \boldsymbol\upnu}(t,n+1)).
	\end{multline*}
\end{Theorem}
In other words,
when $\upnu_n<\upnu_{n+1}$,
the beta polymer admits a Markov swap operator
$p_n^{\textnormal{beta}}$ (in~the sense of Definition~\ref{def:local_transition})
which acts by
splitting the segment
$[Z^{\boldsymbol\upnu}(t,n+1),Z^{\boldsymbol\upnu}(t,n)]\subset[0,1]$
as $1-\tilde B:\tilde B$, and
replacing
$Z^{\boldsymbol\upnu}(t,n)$
by the separating point.
\begin{proof}[Proof of Theorem~\ref{thm:beta_polymer_back}]
	The proof is similar to the case of the $q$-Hahn
	TASEP given in Section~\ref{sub:cond_distr_qhahn}.
	Here we briefly outline the main computations.
	We use the notation~\eqref{eq:n_a_b_notation} which we reproduce here for convenience:
	\begin{gather}
		\label{eq:beta_polymer_proof05}
		\vec n
		=(n_1,\ldots,n_k )=
		(m_1,\ldots,m_\ell,\underbrace{n+1,\ldots,n+1 }_a,\underbrace{n,\ldots,n }_b,m_1',\ldots,m'_{\ell'}).
	\end{gather}
	Applying $p_n^{\textnormal{beta}}$ to the moment formula,
	we will compute moments of
	the form
	\begin{gather}
		\label{eq:beta_polymer_proof1}
		\mathbb{E}^{\mathrm{beta}(\boldsymbol\upnu)}
		\biggl\{
		\big( \tilde B Z(t,n+1)+(1-\tilde B)Z(t,n) \big)^b
		\prod_{\substack{j=1\\n_j \ne n}}^{k}
		Z(t,n_j)\biggr\}.
	\end{gather}
	Expanding
	$( \tilde B Z(t,n+1)+(1-\tilde B)Z(t,n) )^b$
	and using the independence of
	$\tilde B$ from the polymer, we
	have (see, for example, \cite[Lemma 4.1]{CorwinBarraquand2015Beta}
	for the moments of the beta distribution)
	\begin{gather}
		\eqref{eq:beta_polymer_proof1}=
		\sum_{r=0}^b
		\binom{b}{r}
		\frac{(\upnu_n)_r(\upnu_{n+1}\!-\!\upnu_n)_{b-r}}{(\upnu_{n+1})_b}\,
		\mathbb{E}^{\mathrm{beta}(\boldsymbol\upnu)}
		\biggl\{
			Z(t,n)^{r}
			Z(t,n\!+\!1)^{a+b-r}
			\!\!\!\!\!\prod_{\substack{j=1\\n_j \ne n,n+1}}^{k}\!\!\!\!\!
		Z(t,n_j)\biggr\},\!\!\!
		\label{eq:beta_polymer_proof2}
	\end{gather}
	where $(\alpha)_k:=\alpha(\alpha+1)\cdots(\alpha+k-1)$ is the Pochhammer symbol.
	Denote the expectation in~\eqref{eq:beta_polymer_proof2} by
	\begin{gather*}
		g(t;m_1,\ldots,m_\ell,\underbrace{n+1,\ldots,n+1}_{a+b-r},
		\underbrace{n,\ldots,n}_r,m_1',\ldots,m_{\ell'}').
	\end{gather*}
	This expectation
	satisfies the two-body boundary conditions~\eqref{eq:beta_boundary_conditions}
	with the parameter $\upnu=\upnu_{n+1}$.
	Therefore, using the argument from
	\cite[Section 4]{CorwinBarraquand2015Beta}
	(a statement parallel to the $q$-Hahn TASEP's Lemma~\ref{lemma:phi_duality}),
	we can rewrite the sum over $r$ in
	\eqref{eq:beta_polymer_proof2} as the action of the free operators:
	\begin{gather*}
		\prod_{j=1}^b
		\big[\nabla^{\textnormal{beta}}_{\upnu_n/\upnu_{n+1}}\big]_{\ell+a+j}\,
		g(t;m_1,\ldots,m_\ell,\underbrace{n+1,\ldots,n+1}_{a+b},
		m_1',\ldots,m_{\ell'}').
	\end{gather*}
	Finally, each of the
	operators $\nabla^{\textnormal{beta}}_{\upnu_n/\upnu_{n+1}}$ can
	be applied separately under the contour integral in $g(t;\cdot)$ given by
	\eqref{eq:beta_polymer_moments}, and we obtain (with the notation
	$w=z_{\ell+a+j}$, $j=1,\ldots,b$):
	\begin{align*}
		\big[\nabla^{\textnormal{beta}}_{\upnu_n/\upnu_{n+1}}\big]_{\ell+a+j}
		\prod_{i=1}^{n+1}\frac{w}{w-\upnu_i}
		&=	\bigg(\frac{\upnu_n}{\upnu_{n+1}}+
			\bigg( 1-\frac{\upnu_n}{\upnu_{n+1}} \bigg)
			\frac{w}{w-\upnu_{n+1}}\bigg)\prod_{i=1}^{n}\frac{w}{w-\upnu_i}
		\\
&=
		\frac{w-\upnu_{n}}{w-\upnu_{n+1}}\prod_{i=1}^{n}\frac{w}{w-\upnu_i}
	=\prod_{\substack{i=1\\i\ne n}}^{n+1}\frac{w}{w-\upnu_i}.
	\end{align*}
	Thus, we see that~\eqref{eq:beta_polymer_proof1} is equal
	to the expectation
	$\mathbb{E}^{\mathrm{beta}(s_n\boldsymbol\upnu)}( Z(t,n_1)\cdots Z(t,n_k) )$
	with the swapped parameters $\upnu_n\leftrightarrow\upnu_{n+1}$,
	where $(n_1,\ldots,n_k )$ is given by~\eqref{eq:beta_polymer_proof05}.
	Since joint moments determine the distribution of the beta polymer,
	we are done.
\end{proof}

\subsection{Polymer interpretation}

\begin{figure}[htpb]
	\centering
	\includegraphics{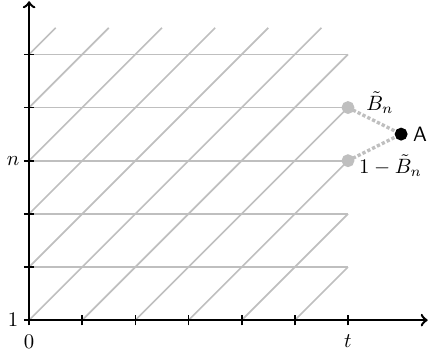}
	\caption{Modified lattice of finite width used to interpret
	$\tilde Z(t,n)$~\eqref{eq:beta_tilde_zn} as a polymer partition function.}
	\label{fig:beta_modified_lattice}
\end{figure}

Let us give a polymer interpretation of Theorem~\ref{thm:beta_polymer_back}
(assuming that $\upnu_n<\upnu_{n+1}$).
First, observe that the quantity
\begin{gather}
	\label{eq:beta_tilde_zn}
\tilde Z(t,n):=\tilde B_n Z(t,n+1)+(1-\tilde B_n)Z(t,n),
\end{gather}
where
$\tilde B_n\sim \mathrm{Beta}(\upnu_{n+1}-\upnu_n,\upnu_n)$, is a
beta polymer type
partition function on a modified lattice.
This modified lattice coincides with the one
in Fig.~\ref{fig:beta_polymer}
in the vertical strip $\{0,1,\ldots,t \}\times \mathbb{Z}_{\ge1}$,
has an additional vertex $\mathsf{A}$, and two additional directed edges
$(t,n)\rightarrow \mathsf{A}$ and $(t,n+1)\to \mathsf{A}$
with weights $1-\tilde B_n$ and $\tilde B_n$, respectively.
The partition function from the line
$\left\{ 0 \right\}\times \mathbb{Z}_{\ge1}$ to $\mathsf{A}$ is precisely
$\tilde Z(t,n)$. See
Fig.~\ref{fig:beta_modified_lattice} for an illustration.

Let us now iterate the swapping in Theorem~\ref{thm:beta_polymer_back}
and
interchange the parameter $\upnu_1$
with~$\upnu_2$, then with $\upnu_3$, and so on up to infinity.
Assume that $\upnu_1<\upnu_2<\cdots $.
Let
$B_n^{(1)}\sim \mathrm{Beta}(\upnu_{n+1}-\upnu_1,\upnu_1)$, $n\in \mathbb{Z}_{\ge1}$,
be independent random variables which are also independent
of the environment $\{B_{t,n}\}$ in the beta polymer.
Define
\begin{gather}
	\label{eq:beta_iterated_map}
	Z^{(1)}(t,n):=B_n^{(1)} Z(t,n+1)+(1-B_n^{(1)})Z(t,n), \qquad n=1,2,\ldots.
\end{gather}
\begin{Proposition}
	\label{prop:beta_iterated_swapping}
	The joint distribution of
	$\{Z^{(1)}(t,n)\}_{n\in \mathbb{Z}_{\ge1}}$
	defined above coincides with
	the joint distribution of
	the beta polymer
	$\{Z^{(\upnu_2,\upnu_3,\ldots )}(t,n)\}_{n\in \mathbb{Z}_{\ge1}}$ at same $t$,
	but
	with the sequence of parameters shifted by one.
\end{Proposition}
\begin{proof}
	Fix $m\in \mathbb{Z}_{\ge1}$.
	Observe that the
	joint distribution of
	$(Z^{(\upnu_1,\upnu_2,\ldots )}(t,m),Z^{(\upnu_1,\upnu_2,\ldots )}(t,m+1))$
	does not depend on the order of the parameters
	$\upnu_1,\ldots,\upnu_m $.
	Therefore, applying~\eqref{eq:beta_iterated_map}
	with $n=m$ and using Theorem~\ref{thm:beta_polymer_back}
	makes the new random variable
	$Z^{(1)}(t,m)$ a beta polymer partition function
	with parameters $(\upnu_2,\ldots,\upnu_n,\upnu_{n+1} )$.
	The statement about joint distributions is
	obtained by sequential application of
	this argument
	for $m=1,2,\ldots $.
\end{proof}

\begin{figure}[htpb]
	\centering
	\includegraphics{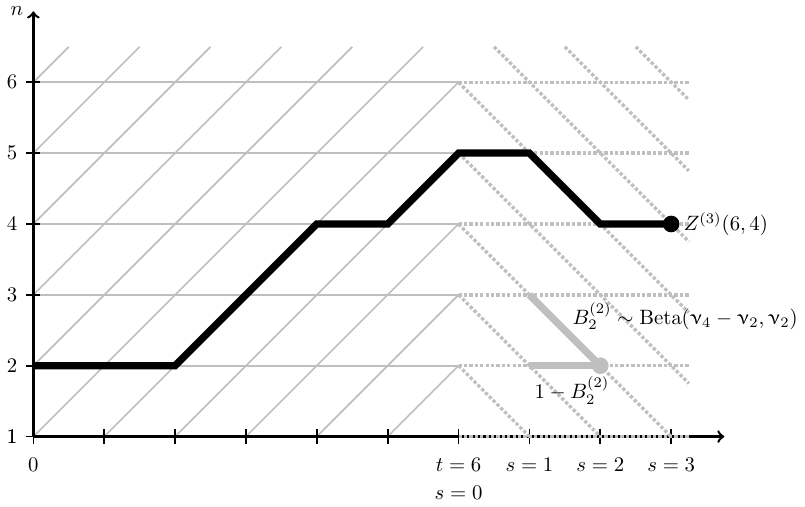}
	\caption{The lattice used to define the beta polymer partition functions
	$Z^{(s)}(t,n)$ with shifted parameter sequences.}
	\label{fig:beta_modified_lattice2}
\end{figure}

The quantities $\{Z^{(1)}(t,n)\}_{n\in \mathbb{Z}_{\ge1}}$
can be interpreted as beta polymer type partition functions, too.
Moreover, let us further iterate
Proposition~\ref{prop:beta_iterated_swapping},
and introduce independent random variables
\begin{gather*}
	B^{(s)}_n\sim \mathrm{Beta}	(\upnu_{n+s}-\upnu_s,\upnu_s),\qquad
	s=1,2,\ldots.
\end{gather*}
Define $Z^{(s)}(t,n)$ to be the polymer partition
function from the line $\left\{ 0 \right\}\times \mathbb{Z}_{\ge1}$
to the point $(s+t,n)$, $s\in \mathbb{Z}_{\ge1}$,
in the modified strict-weak lattice
which coincides with the original lattice
in Fig.~\ref{fig:beta_polymer} in the vertical strip
$\left\{ 0,1,\ldots,t \right\}\times \mathbb{Z}_{\ge1}$.
To the right of this strip,
the modified lattice is
made out of down-right diagonal and horizontal
edges with the weights
$B_n^{(s)}$ on each
$(t+s-1,n+1)\to(t+s,n)$, and
$1-B_n^{(s)}$ on each $(t+s-1,n)\to(t+s,n)$.
See Fig.~\ref{fig:beta_modified_lattice2} for~an~illustration.

{\sloppy\begin{Proposition}
	\label{prop:joint_shifted}
	For any fixed
	$s$ and $t$,
	the joint distribution of
	the partition functions
	$\{Z^{(s)}(t,n)\}_{n\in \mathbb{Z}_{\ge1}}$ on the modified lattice
	coincides with the joint distribution of the
	beta polymer partition functions
	$\{Z^{(\upnu_{s+1},\upnu_{s+2},\ldots )}(t,n)\}_{n\in \mathbb{Z}_{\ge1}}$
	with the same $t$,
	but with the parameter sequence~$\boldsymbol\upnu$
	shifted by $s$.
\end{Proposition}

}

\subsection{Zero-temperature limit}

Under a limit transition when the parameters
of the beta random variables go to zero,
the beta polymer model
turns into a first passage percolation type model.
First, we recall the scaling:
\begin{Lemma}[{\cite[Lemma 5.1]{CorwinBarraquand2015Beta}}]
	Let $\alpha,\beta>0$, and $B_\varepsilon\sim\mathrm{Beta(\varepsilon \alpha,
	\varepsilon\beta)}$.
	Then, as $\varepsilon\searrow0$, we have convergence in distribution:
	\begin{gather*}
		\left( -\varepsilon\log B_\varepsilon,
		-\varepsilon\log(1-B_\varepsilon)\right)
		\to
		(\xi E_\alpha,(1-\xi)E_\beta).
	\end{gather*}
	Here $\xi\in\left\{ 0,1 \right\}$ is the Bernoulli random variable with
	$\mathbb{P}(\xi=1)=\frac{\beta}{\alpha+\beta}$,
	and $(E_\alpha,E_\beta)$ are exponential random variables with parameters
	$\alpha$ and $\beta$ $($that is, means $\alpha^{-1}$ and $\beta^{-1})$ which are independent of $\xi$.
\end{Lemma}

\begin{figure}[htpb]
	\centering
	\includegraphics{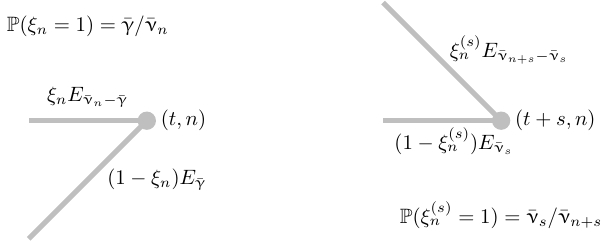}
	\caption{Edge weights $t_e$ in the zero-temperature limit.
	The $\xi$'s are independent Bernoulli random variables with given parameters,
	and all the $E_\alpha$'s are exponential random variables independent
	of the Bernoulli ones.}
	\label{fig:exp_weights}
\end{figure}

We will take the scaling limit of the beta polymer model
as $\upnu_n=\varepsilon\bar\upnu_n$,
$\upgamma=\varepsilon\bar \upgamma$,
where $\bar \upnu_n>\bar\upgamma>0$ for all
$n$, and $0<\bar\upnu_1<\bar\upnu_2<\cdots $.
The edge weights in the lattice in
Fig.~\ref{fig:beta_modified_lattice2}
turn into the ones given in Fig.~\ref{fig:exp_weights}.

Denote by $F^{(s)}(t,n)$ the first-passage time
from the line $\left\{ 0 \right\}\times \mathbb{Z}_{\ge1}$ to the
point $(s+t,n)$ in~the modified lattice
\begin{gather*}
	F^{(s)}(t,n):=\min_{\pi\colon\left\{ 0 \right\}\times \mathbb{Z}_{\ge1}\to (t+s,n)}
	\sum_{e\in \pi} t_e,
\end{gather*}
where the directed paths $\pi$ are as in Fig.~\ref{fig:beta_modified_lattice2},
and the edge weights are given in Fig.~\ref{fig:exp_weights}.
If $s=0$, then we mean the unmodified first-passage time
(as studied in~\cite{CorwinBarraquand2015Beta}).
Above the main diagonal (i.e., for $n>t$) we have $F^{(0)}(t,n)=0$
because due to the presence of the Bernoulli components,
there always exists a path with zero total weight
between $(t,n)$, $n>t$, and the vertical axis $\{0 \}\times \mathbb{Z}_{\ge1}$.

For the first-passage percolation model, an analogue of Proposition~\ref{prop:joint_shifted}
holds:
\begin{Proposition}	
	For fixed $s$, $t$,
	the joint distribution of the first-passage times
	$\{F^{(s)}(t,n)\}_{n\in \mathbb{Z}_{\ge1}}$
	with parameters $\bar\upnu_1<\bar\upnu_2<\cdots $
	in the modified lattice
	is the same as that of the unmodified ones
	$\big\{F^{(0)}(t,n)\big\}_{n\in \mathbb{Z}_{\ge1}}$,
	but with the shifted sequence of parameters
	$\bar\upnu_{s+1}<\bar\upnu_{s+2}<\cdots $.
\end{Proposition}

\subsection*{Acknowledgements}

I am grateful to Vadim Gorin for helpful discussions,
and to Matteo Mucciconi and Axel Saenz for remarks on the first version of the manuscript.
I am grateful to the organizers of the workshop
``Dimers, Ising Model, and their Interactions''
and the support of the Banff International
Research Station where a part of this work was done.
The work was partially supported by the NSF grant DMS-1664617.

\pdfbookmark[1]{References}{ref}
\LastPageEnding

\end{document}